\documentclass[leqno]{article}
\usepackage{amssymb,amsfonts,amsmath,theorem,pstricks,graphicx,subfig,longtable,float}

\newcommand{\Tor}{\text{Tor}}

\newtheorem{lemma}{Lemma}

\newtheorem{theorem}{Theorem}
\newtheorem{conjecture}{Conjecture}
{\theorembodyfont{\rmfamily}}

\begin{document}

\title{Newton polygons and curve gonalities}
\author{Wouter Castryck\footnote{Supported by the Fund for Scientific Research Flanders (F.W.O.\ Vlaanderen) }
\ and Filip Cools${}^*$}

\maketitle

\begin{abstract}
\noindent We give a combinatorial upper bound for the gonality of a
curve that is defined by a bivariate Laurent polynomial with given Newton polygon. We
conjecture that this bound is generically attained, and provide proofs in a considerable number of
special cases. One proof technique uses recent work of M.\ Baker on linear systems on graphs, by means of which we reduce our conjecture
to a purely combinatorial statement.\\

\end{abstract}

\section{Introduction}

The most renowned birational invariant of an algebraic curve over $\mathbb{C}$ is its geometric genus.
Although it enjoys the plastic description as the number of handles on the
corresponding Riemann surface, some high-tech machinery is needed to give a rigorous definition.
E.g., one nowadays approach is to define the geometric genus as the $\mathbb{C}$-dimension of the Riemann-Roch space
associated to a canonical divisor $K_C$ on $C$.

On the other hand, at the end of the $19$th century already, H.\ Baker shared the following elementary observation.
Let $f \in \mathbb{C}[x^{\pm 1}, y^{\pm 1}]$ be an irreducible Laurent polynomial defining
a curve $U(f) \subset \mathbb{T}^2$, where $\mathbb{T}^2 = \left( \mathbb{C} \setminus \{0\} \right)^2$ is
the two-dimensional torus over $\mathbb{C}$.
Let $\Delta(f)$ be the Newton polygon of $f$. It is an instance of
a \emph{lattice polygon}, by which we mean the convex hull in $\mathbb{R}^2$ of a finite subset of $\mathbb{Z}^2$.
The \emph{dimension} of a lattice polygon $\Delta$ is the minimal dimension of an affine subspace of $\mathbb{R}^2$ containing
$\Delta$. By the \emph{interior} of $\Delta$ we mean the topological interior if $\Delta$ is two-dimensional, and the empty
set if it is of strictly lower dimension. Points of $\mathbb{Z}^2$ will be called \emph{lattice points}. Then:
\begin{theorem}[Baker, 1893] \label{bakerinequality}
The geometric genus of $U(f)$ is at most the number of lattice points in the interior of $\Delta(f)$.
\end{theorem}
\noindent \textsc{Proof.} This can be found in \cite{HBaker}. See \cite{Beelen} for a more modern proof. \hfill $\blacksquare$\\

Generically, Baker's bound is sharp.

\begin{theorem}[Khovanski\u{\i}, 1977] \label{Khovanskii}
Let $\Delta \subset \mathbb{R}^2$ be a two-dimensional lattice polygon. The set of irreducible Laurent polynomials
$f \in \mathbb{C}[x^{\pm 1}, y^{\pm 1}]$ for which $\Delta(f) = \Delta$ and
the bound in Theorem~\ref{bakerinequality} is attained,
is Zariski dense in the space of Laurent polynomials $f \in \mathbb{C}[x^{\pm 1}, y^{\pm 1}]$ for which $\Delta(f) \subset \Delta$.
\end{theorem}

\noindent \textsc{Proof.} See \cite{Khovanskii}. Khovanski\u{\i} actually proved
something much stronger, which we will state in Section~\ref{canonicalimage}. \hfill $\blacksquare$\\

Because of all this, one defines the \emph{genus} of a lattice polygon as the number of
lattice points in its interior.

The near-miraculous appearance of the interior lattice points of $\Delta(f)$ was
secularized with
the advent of tropical geometry. Loosely stated, by subdividing $\Delta(f)$ into triangles of area $1/2$ and taking the dual
spine, one obtains a graph whose handles are in one-to-one correspondence with the interior lattice points of $\Delta(f)$.
By considering the graph as a piece-wise linear limit of $U(f)$ in an appropriate way, one realizes that one has actually visualized
the handles of the Riemann surface. This construction will be reviewed somehow in Sections~\ref{toricgeometry} and \ref{graphgonalities} below.

The aim of this article is to give analogues of Theorems~\ref{bakerinequality} and \ref{Khovanskii} for what is, arguably, the second most renowned
birational invariant of an algebraic curve: its gonality.
Our results will be partly conjectural.
To state them, we need the following terminology.
A \emph{$\mathbb{Z}$-affine transformation}
is a map $\mathbb{R}^2 \rightarrow \mathbb{R}^2 : x \mapsto Ax + b$ with
$A \in \text{GL}_2(\mathbb{Z})$ and $b \in \mathbb{Z}^2$.
Two lattice polygons $\Delta, \Delta'$ are called \emph{equivalent} if there is
a $\mathbb{Z}$-affine transformation $\varphi$ such that $\varphi(\Delta) = \Delta'$ (notation: $\Delta \cong \Delta'$).
The \emph{lattice width} of a non-empty lattice polygon $\Delta$ is
the smallest integer $s \geq 0$ such
that there is a $\mathbb{Z}$-affine transformation $\varphi$ for which $\varphi(\Delta)$ is
contained in the horizontal strip
\[ H_0^s = \left\{ \left. \, (x,y) \in \mathbb{R}^2 \, \right| \, 0 \leq y \leq s \, \right\}.\]
It will be denoted by $\text{lw}(\Delta)$. It is convenient to define $\text{lw}(\emptyset) = -1$.
If $\Delta$ is a lattice polygon, then for any integer $d \geq 0$, the polygon $d \Delta$ denotes the corresponding Minkowski multiple.
We will denote the standard $2$-simplex in $\mathbb{R}^2$ by $\Sigma$. Thus,
$d\Sigma$ is the Newton polygon of a generic degree $d$ polynomial. We use $\Upsilon$ to denote
$\text{Conv} \{ (-1,-1),(1,0),(0,1) \}$.

Then our analogues of Theorems~\ref{bakerinequality} and \ref{Khovanskii} read:

\begin{theorem} \label{goninequality}
The gonality of $U(f)$ is at most $\emph{lw}(\Delta(f))$. If
$\Delta(f)$ is equivalent to $d\Sigma$
for some $d \geq 2$, or
to $2 \Upsilon$,
then it is at most $\emph{lw}(\Delta(f)) - 1$.
\end{theorem}

\begin{conjecture} \label{gongeneric}
Let $\Delta \subset \mathbb{R}^2$ be a two-dimensional lattice polygon. The set of irreducible Laurent polynomials
$f \in \mathbb{C}[x^{\pm 1}, y^{\pm 1}]$ for which $\Delta(f) = \Delta$ and
the (sharpest applicable) bound in
Theorem~\ref{goninequality} is attained, is Zariski dense in the space of Laurent polynomials $f \in \mathbb{C}[x^{\pm 1}, y^{\pm 1}]$ for which $\Delta(f) \subset \Delta$.
\end{conjecture}

The article is organized as follows.

In Section~\ref{upperbound}, we prove Theorem~\ref{goninequality}.

In Section~\ref{interior}, we give a reformulation of Conjecture~\ref{gongeneric} that focuses
on the convex hull of the interior lattice points of $\Delta$, rather than on $\Delta$ itself. In doing so, the polygons $d\Sigma$ become
ruled out as special instances.

In Section~\ref{toricsurfaces}, we review how to associate
a toric surface $\text{Tor}(\Delta)$ to a lattice polygon $\Delta$, and how, in general, $\text{Tor}(\Delta(f))$ naturally appears
as an ambient space for the complete non-singular model of $U(f)$.

In Section~\ref{canonicalimage}, we prove
Conjecture~\ref{gongeneric} for all $\Delta$ for which $\text{lw}(\Delta) \leq 4$ (including $\Delta \cong 2\Upsilon$), by analyzing
the canonical image of $U(f)$. We briefly report on a
computer experiment supporting Conjecture~\ref{gongeneric}
for all lattice polygons up to genus $13$, thereby relying on Green's canonical conjecture.


In Section~\ref{previousresults}, we see how previous results by Kawaguchi, Martens, and Namba
prove Conjecture~\ref{gongeneric} in a considerable number of additional cases.

In Section~\ref{toricgeometry}, we review the process of degenerating a toric surface $\text{Tor}(\Delta)$
according to a regular subdivision of $\Delta$, and use this to deform a sufficiently generic $U(f)$ along
with $\text{Tor}(\Delta)$ into a union of irreducible curves. As a by-product, we obtain a vast class of examples
of strongly semi-stable arithmetic surfaces.

In Section~\ref{graphgonalities}, we encode
the combinatorial configuration of this union of irreducible curves in a graph, and we apply
recent results due to M.\ Baker \cite{MBaker} to obtain a lower bound for
the gonality of $U(f)$.

In Section~\ref{purelycombinatorial} we conjecture that, in this way, one can always meet the
upper bound of Theorem~\ref{goninequality}. This reduces
Conjecture~\ref{gongeneric} to a purely combinatorial (albeit a priori stronger) statement.
We prove this statement (and hence Conjecture~\ref{gongeneric}) for an interesting class of
lattice polygons, thereby partly confirming and partly extending the results
of Sections~\ref{canonicalimage} and~\ref{previousresults}.\\

\noindent \emph{Acknowledgements.} We are very
grateful to an anonymous referee for his valuable comments, that led to the proof
of Theorem~\ref{theoremnewfamily}.
We would also like to thank Marc Coppens, Hendrik Hubrechts, Bjorn Poonen, Jan Schepers, Jan Tuitman and Wim Veys for
some helpful discussions. The first author thanks the Massachusetts Institute of Technology for
its hospitality.

\section{The lattice width as an upper bound} \label{upperbound}

In this section, we prove Theorem~\ref{goninequality}.\\

\noindent \textsc{Proof of Theorem~\ref{goninequality}.} It is clear that if $\Delta(f)$
is contained in a horizontal strip
\[ H_0^s = \left\{ \left. \, (x,y) \in \mathbb{R}^2 \, \right| \, 0 \leq y \leq s \, \right\},\]
then the rational map
\[ U(f) \rightarrow \mathbb{A}^1 : (x,y) \mapsto x \]
is of degree at most $s$. Now every $\mathbb{Z}$-affine transformation $\varphi : \mathbb{R}^2 \rightarrow \mathbb{R}^2$ acts
on $f$ as follows: if
\[ f = \sum_{(i,j) \in \Delta \cap \mathbb{Z}^2} c_{ij} (x,y)^{(i,j)}, \quad \text{then } f^\varphi = \sum_{(i,j) \in \Delta \cap \mathbb{Z}^2} c_{ij} (x,y)^{\varphi(i,j)} \]
(where we use multi-index notation). It is clear that $\Delta(f^\varphi) = \varphi(\Delta(f))$ and that $U(f) \cong U(f^\varphi)$.
The upper bound $\text{lw}(\Delta)$ follows immediately.

Now suppose that $\Delta(f) \cong d\Sigma$ for some integer $d \geq 2$. Hence we can
assume that $f \in \mathbb{C}[x,y]$ is a dense degree $d$ polynomial, whose homogenization $F$ with
respect to a new variable $z$ defines a curve $V(F)$ in $\mathbb{P}^2 = \text{Proj} \, \mathbb{C}[x,y,z]$.
A projective transformation takes us to a curve $V(F')$ containing the point $(0:1:0)$. Dehomogenizing $F'$ with respect
to $z$ yields a polynomial $f' \in \mathbb{C}[x,y]$ whose Newton polygon is contained in
\[ \text{Conv} \{ (0,0),(d,0),(1,d-1),(0,d-1) \}.\]
Thus, $\Delta(f')$ is of lattice width at most $d-1$. Hence the gonality of $U(f)$ is at most $d-1$. On the other hand, $\text{lw}(d\Sigma) = d$ for
all integers $d \geq 0$. Indeed, clearly
$\text{lw}(d\Sigma) \leq d$, and equality follows from the fact that each edge of $d\Sigma$ contains $d+1$ lattice points.

Finally, suppose that $\Delta(f) \cong 2 \Upsilon$. By Theorem~\ref{bakerinequality}, $U(f)$
has geometric genus at most $4$, and it is classical that this implies the gonality to be at most $3$, see e.g.\ \cite{KleimanLaksov}.
On the other hand, $\text{lw}(2 \Upsilon) = 4$. Indeed, clearly $\text{lw}(2 \Upsilon) \leq 4$, and equality follows
from the fact that the convex hull of the interior lattice points of $2 \Upsilon$ has itself an interior lattice point.
But this would be impossible if $\text{lw}(2 \Upsilon) \leq 3$. \hfill $\blacksquare$

\section{The interior lattice polygon} \label{interior}

Our two exceptional cases $d \Sigma (d \geq 2)$ and $2 \Upsilon$ are of a very different kind.
In the first case, one is able clip off a vertex in such a way that it reduces the lattice width,
without affecting the geometry of $U(f)$. For $2 \Upsilon$, such a trick is impossible, since clipping off
a vertex would necessarily mean reducing the number of interior lattice points. Hence this would affect the
generic genus of $U(f)$.

In this section we will deduce an equivalent formulation of Conjecture~\ref{gongeneric}, in which
the polygons $d\Sigma$ are no longer exceptional cases. This is done by focusing on the \emph{interior lattice polygon}, rather
than on the polygon itself. For any lattice polygon $\Delta$, the interior lattice polygon
$\Delta^{(1)}$ is defined as the convex hull of the interior lattice points of $\Delta$. Somehow dually, one can
consider the \emph{relaxed polygon}. That is, let $\Delta$ be a two-dimensional lattice polygon, and write it
as a finite intersection of half-planes
\[ \Delta = \bigcap_i \left\{ (x,y) \in \mathbb{R}^2 \, | \, a_ix + b_iy \leq c_i \right\},\]
where $a_i,b_i,c_i \in \mathbb{Z}$ and $\gcd(a_i,b_i) = 1$ for all $i$. Then the relaxed polygon is defined as
\[ \Delta^{(-1)} = \bigcap_i \left\{ (x,y) \in \mathbb{R}^2 \, | \, a_ix + b_iy \leq c_i + 1 \right\}.\]
Not every lattice polygon can be written as $\Delta^{(1)}$ for some larger lattice polygon $\Delta$.
Also, if $\Delta$ is a two-dimensional lattice polygon $\Delta$, then $\Delta^{(-1)}$ need not
be a lattice polygon: it may take vertices outside the lattice.
The following statement connects and controls both phenomena.

\begin{lemma} \label{haaseschicho}
Let $\Delta$ be a two-dimensional lattice polygon. Then $\Delta = \Gamma^{(1)}$
for a lattice polygon $\Gamma$ if and only if $\Delta^{(-1)}$ is a lattice polygon. Moreover,
if $\Delta^{(-1)}$ is a lattice polygon, then it is maximal
(with respect to inclusion) among all lattice polygons $\Gamma$ for which $\Gamma^{(1)} = \Delta$.
\end{lemma}

\noindent \textsc{Proof.} This is due to \cite[Section~2.2]{Koelman}.
Recently, this was rediscovered by Haase and Schicho \cite[Lemmata~9~\&~11]{HaaseSchicho}. \hfill $\blacksquare$\\

The main result of this section is the following relationship between $\text{lw}(\Delta)$ and $\text{lw}(\Delta^{(1)})$.
It was discovered independently (and almost simultaneously) by Lubbes and Schicho \cite[Theorem~13]{LubbesSchicho}.

\begin{theorem} \label{lwDelta1} Let $\Delta$ be a two-dimensional lattice polygon.
Then
\[ \emph{lw}(\Delta) = \emph{lw}(\Delta^{(1)}) + 2,\]
unless $\Delta \cong d \Sigma$ for some integer $d \geq 2$, in which
case $\emph{lw}(\Delta) = \emph{lw}(\Delta^{(1)}) + 3 = d$.
\end{theorem}

\noindent \textsc{Proof.}
First, it is clear that $\Delta^{(1)}$ can be caught in a horizontal strip of
width $\text{lw}(\Delta) - 2$, from which
\begin{equation} \label{lwineq}
  \text{lw}(\Delta^{(1)}) \leq \text{lw}(\Delta) - 2.
\end{equation}
Second, in Section~\ref{upperbound} we saw that $\text{lw}(d\Sigma) = d$ for all integers $d \geq 0$. Third, we
have that $(d \Sigma)^{(1)} \cong (d-3)\Sigma$ for all integers $d \geq 3$, and that
$(2\Sigma)^{(1)} = \emptyset$. So $\text{lw}(d \Sigma) = \text{lw}(d \Sigma^{(1)}) + 3 = d$ for all
$d \geq 2$.
We therefore conclude that it suffices to prove: if the inequality in (\ref{lwineq}) is
strict, then $\Delta \cong d\Sigma$ for some integer $d \geq 2$.

For technical reasons, we first get rid of the following cases.
\begin{itemize}
\item $\text{lw}(\Delta^{(1)})=-1$, i.e.\ $\Delta$ contains no interior lattice points.
Then either $\Delta \cong 2\Sigma$ or $\Delta$ is a so-called Lawrence prism, see \cite[Thm 4.1.2]{Koelman} or the generalized statement of \cite[Thm 2.5]{BatyrevNill}. Since Lawrence prisms have lattice width $1$, the result follows.
\item $\text{lw}(\Delta^{(1)})=0$. Then $\Delta$ is a so-called elliptic or hyperelliptic
lattice polygon. These have been classified in \cite[Thm 4.2.3 and Sec 4.3]{Koelman}, from
which it follows that either $\Delta \cong 3\Sigma$, or $\text{lw}(\Delta) = 2$.
\item $\text{lw}(\Delta^{(1)})=1$. Then $\Delta^{(1)}$ must be a Lawrence prism,
and using Lemma~\ref{haaseschicho} one concludes that either $\Delta \cong 4\Sigma$, or
$\text{lw}(\Delta) = 3$.
\item $\Delta^{(1)}$ can be caught in a $3$-by-$3$ lattice square.
These cases can be exhaustively verified using Lemma~\ref{haaseschicho}.
\end{itemize}

To deal with the general case, we apply a $\mathbb{Z}$-affine transformation to catch
$\Delta^{(1)}$ in
the horizontal strip
\[ H = H_0^{\text{lw}(\Delta^{(1)})} = \left\{ \left. (x,y) \in \mathbb{R}^2 \, \right| \, 0 \leq y \leq \text{lw}(\Delta^{(1)}) \, \right\}. \]
Since we assume (\ref{lwineq}) strict, $\Delta$ must then contain at least one vertex $v$ outside the strip
\[ \left\{ \left. (x,y) \in \mathbb{R}^2 \, \right| \, -1 \leq y \leq \text{lw}(\Delta^{(1)}) + 1 \, \right\}. \]
We may assume that $v = (0,-k)$ with $k \geq 2$.

We will first prove that $k=2$. Since $\Delta$ contains no interior
lattice points on the line $y=-1$, it must intersect this line
inside an interval $[\alpha, \alpha+1]$ for some $\alpha \in \mathbb{Z}$.
Let $\sigma$ be the cone with top $v$, whose
rays pass through $(\alpha,-1)$ and $(\alpha+1,-1)$ respectively.
Although $\Delta$ need (a priori) not be contained in $\sigma$,
the part of $\Delta$ that lies on or above the line $y=-1$ must be. In particular,
$\Delta^{(1)}$ will be contained in the open cone $\sigma^\circ$.
Modulo horizontally skewing and flipping if necessary, we can assume that
\begin{equation} \label{boundalpha}
  0 \leq \alpha \leq \left\lfloor \frac{k-2}{2} \right\rfloor.
\end{equation}
Then $\Delta^{(1)} \subset \sigma^\circ \cap H$ is contained in the interior of
the vertical strip
\[ V = \left\{  (x,y) \in \mathbb{R}^2 \, \left| \, \frac{\alpha k}{k-1} \leq x \leq \frac{(\alpha + 1)(\text{lw}(\Delta^{(1)})+k)}{k-1} \, \right. \right\} \]
which has width
\[ \frac{(\alpha +1)}{k-1} \text{lw}(\Delta^{(1)}) + \frac{k}{k-1} \leq \frac{ \left\lfloor \frac{k}{2} \right\rfloor}{k-1} \text{lw}(\Delta^{(1)}) + \frac{k}{k-1}.\]
By definition of the lattice width
\[ \text{lw}(\Delta^{(1)}) < \frac{ \left\lfloor \frac{k}{2} \right\rfloor}{k-1} \text{lw}(\Delta^{(1)}) + \frac{k}{k-1}.\]

\begin{figure}[ht]
\centering
\includegraphics[height=4cm]{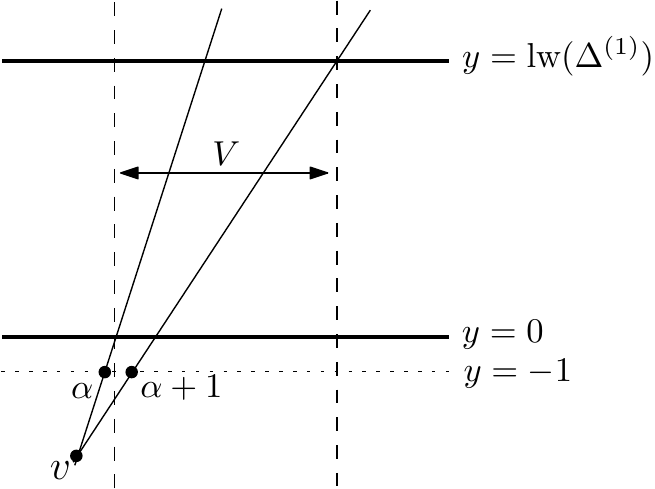}
\end{figure}

\noindent This is impossible for $k \geq 3$ as soon as $\text{lw}(\Delta^{(1)}) \geq 4$.
If $\text{lw}(\Delta^{(1)}) \in \{2, 3 \}$, we find that
$k \geq 3$ would cause
$\Delta^{(1)} \subset H \cap V$ to be caught in a $3$-by-$3$ square, a
case covered in the list above.

Next, note that $v$ is the only vertex of $\Delta$ on the line $y=-2$. Indeed, along
with a vertex of $\Delta^{(1)}$ on the line $y=0$, two such vertices
would span a triangle which by Pick's theorem would need
to contain a lattice point on the line $y=-1$. But this would
be an interior lattice point of $\Delta$: a contradiction. As a consequence,
$v$ is the only lattice point of $\Delta$ lying strictly under the line $y=-1$,
hence $\Delta \subset \sigma$.
By (\ref{boundalpha}),
we can assume that $\sigma$ has top $(0,-2)$, and that
its rays pass through $(0,-1)$ and $(1,-1)$, respectively.

The next step is to prove that $\Delta$ cannot have a vertex lying strictly above the line
$y = \text{lw}(\Delta^{(1)}) + 1$. Suppose that there is such a vertex $w$. By symmetry
of arguments, it must be unique and lying on the line $y = \text{lw}(\Delta^{(1)}) + 2$.
Since $w$ must be contained in $\sigma$, its $x$-coordinate
must be among $0, \dots, \text{lw}(\Delta^{(1)}) + 4$.
As before, there must exist
an integer $\beta$ such that
$\Delta$ is contained in the cone $\tau$ with top $w$, whose
rays pass through $(\beta, \text{lw}(\Delta^{(1)}) + 1)$ and $(\beta + 1, \text{lw}(\Delta^{(1)}) + 1)$.
Using that $\tau$ must contain $v$, a case-by-case analysis show that
$\sigma \cap \tau$ is too small for the lattice width
of $\Delta \subset \sigma \cap \tau$ to exceed $\text{lw}(\Delta^{(1)}) + 2$, a
contradiction with the assumed strict inequality in (\ref{lwineq}).

Overall, we obtained that $\Delta$
must be contained in the triangle spanned by $(0,-2)$, $(0,\text{lw}(\Delta^{(1)}) + 1)$
and $(\text{lw}(\Delta^{(1)}) + 3,\text{lw}(\Delta^{(1)}) + 1)$,
which is a copy of $((\text{lw}(\Delta^{(1)}) + 3)\Sigma$. If any of these three
points does not appear as a vertex, $\Delta$ can be seen
to have lattice width at most $\text{lw}(\Delta^{(1)}) + 2$, a contradiction
with the strict inequality in (\ref{lwineq}). Hence $\Delta$ must be the full triangle, QED. \hfill $\blacksquare$\\

Note that Theorem~\ref{lwDelta1} yields an algorithm for recursively computing $\text{lw}(\Delta)$.

We can now rephrase Conjecture~\ref{gongeneric} as follows.\\

\noindent \textbf{Conjecture~\ref{gongeneric} (equivalent formulation)} \emph{Let $\Delta \subset \mathbb{R}^2$ be a
two-dimensional lattice polygon,
and let $S \subset \mathbb{C}[x^{\pm 1}, y^{\pm 1}]$ be the set of irreducible Laurent polynomials for which $\Delta(f) = \Delta$
and
\begin{itemize}
  \item $U(f)$ has gonality $\emph{lw}(\Delta^{(1)}) + 2$ if $\Delta \not \cong 2\Upsilon$,
  \item $U(f)$ has gonality $3$ if $\Delta \cong 2\Upsilon$.
\end{itemize}
Then $S$ is Zariski dense in the space of Laurent polynomials $f \in \mathbb{C}[x^{\pm 1},y^{\pm 1}]$ for which $\Delta(f) \subset \Delta$.}\\

\noindent \textsc{Proof of equivalence.} This follows directly from Theorem~\ref{lwDelta1}. \hfill $\blacksquare$

\section{Toric surfaces as ambient spaces} \label{toricsurfaces}

We give a brief, notation-fixing overview of the geometry of toric surfaces. In Section~\ref{toricgeometry},
the material below will be put in a bigger framework.

Let $\Delta \subset \mathbb{R}^2$ be a two-dimensional lattice polygon. Let $S$ be
the set of lattice points of $\Delta$. Then we have an injective morphism
\[ \phi : \mathbb{T}^2 \hookrightarrow \mathbb{P}^{|S| - 1} : (x,y) \mapsto \left( x^iy^j \right)_{(i,j) \in S}.\]
The Zariski closure of the image is by definition the \emph{toric surface} $\text{Tor}(\Delta)$. If
$X_{i,j}$ denotes the projective coordinate of $\mathbb{P}^{|S| - 1}$ corresponding to
$(i,j) \in S$, then all binomials of the form
\[ \prod_{k=1}^n X_{i_k,j_k} - \prod_{k=1}^n X_{i'_k,j'_k} \]
for which
\[ \sum_{k=1}^n (i_k,j_k) = \sum_{k=1}^n (i'_k,j'_k) \]
are zero on $\text{Tor}(\Delta)$, and in fact these generate the homogeneous ideal of $\text{Tor}(\Delta)$. In practice, it suffices
to consider relations of degree $n \leq 3$, and even $n \leq 2$ if $\# (\partial \Delta \cap \mathbb{Z}^2) > 3$
by a result of Koelman \cite{Koelman2}.

The faces $\tau \subset \Delta$ (vertices, edges, and $\Delta$ itself) naturally partition $\text{Tor}(\Delta)$
into sets
\[ O(\tau) = \left\{ \, \left. (\alpha_{i,j})_{(i,j) \in S} \in \text{Tor}(\Delta) \, \right| \, \alpha_{i,j} \neq 0 \text{ if and only if }
(i,j) \in \tau \right\},\]
which are called the \emph{toric orbits}. Note that $O(\Delta) = \phi(\mathbb{T}^2)$. More generally, one has $O(\tau) \cong \mathbb{T}^{\dim \tau}$.
One can show that $\text{Tor}(\Delta)$ is non-singular, except possibly at the zero-dimensional toric orbits.

Write $f = \sum_{(i,j) \in S} c_{i,j} x^iy^j$. Then $\phi(U(f))$ satisfies
\[ \sum_{(i,j) \in S} c_{i,j} X_{i,j},\]
so it embeds into a hyperplane section of $\text{Tor}(\Delta)$.
More generally, if $\Delta = d \Delta'$ for an integer $d \geq 1$
and a lattice polygon $\Delta'$, then $U(f)$ can be embedded
in a degree $d$ hypersurface section of $\text{Tor}(\Delta')$. Generically, this hyperplane/hypersurface section will
be a complete non-singular model of $U(f)$.
A sufficient condition is that $f$ is \emph{non-degenerate with
respect to its Newton polygon} $\Delta$, meaning that for each face $\tau \subset \Delta(f)$ (vertices, edges, and $\Delta(f)$ itself),
the system
\[ f_\tau = x \frac{\partial f_\tau}{\partial x} = y \frac{\partial f_\tau}{\partial y} = 0\]
has no solutions in $\mathbb{T}^2$. Here, $f_\tau$ is obtained from $f$ by only considering those terms whose exponent vector is contained in $\tau$.
Geometrically, non-degeneracy can be rephrased as follows: the Zariski closure of $\phi(U(f))$ has no singular points in $O(\Delta)$, intersects
the one-dimensional toric orbits transversally, and does not contain the zero-dimensional toric orbits.
This is indeed a generic condition, since non-degeneracy can be rephrased in terms of the non-vanishing of
a certain integral polynomial expression in the coefficients $c_{i,j}$, realized as a product of principal $A$-discriminants in the sense of \cite{GKZ}.
See \cite[Section~2]{CastryckVoight} for some additional details.

\section{Polygons of small lattice width} \label{canonicalimage}

In this section we prove Conjecture~\ref{gongeneric} for lattice polygons $\Delta$ satisfying $\text{lw}(\Delta) \leq 4$.
By the results of Section~\ref{interior}, it suffices to do this for two-dimensional lattice polygons $\Delta$ for which $\text{lw}(\Delta^{(1)}) \leq 2$.

Conjecture~\ref{gongeneric} is automatic in case $\text{lw}(\Delta^{(1)}) = -1$, since by definition, genus $0$ curves have gonality $1$.
Next, Theorem~\ref{lwDelta1} immediately implies Conjecture~\ref{gongeneric} for polygons $\Delta$ for which
$\text{lw}(\Delta^{(1)}) = 0$.
Indeed, by Theorem~\ref{Khovanskii}, our curve $U(f)$ will generically have genus at least $1$, hence
gonality at least $2$.
But by Theorem~\ref{lwDelta1}, either $\text{lw}(\Delta) = 2$ or $\Delta \cong 3\Sigma$,
and the statement follows from Theorem~\ref{goninequality}.

In order to extend this to a proof for the cases where
$\text{lw}(\Delta^{(1)}) \in \{-1,0,1,2\}$, including $\Delta \cong 2 \Upsilon$, we
need the following refined version of Theorem~\ref{Khovanskii}.\\

\noindent \textbf{Theorem~\ref{Khovanskii} (Khovanski\u{\i}, 1977, refined formulation)} \emph{
Let $f \in \mathbb{C}[x^{\pm 1}, y^{\pm 1}]$ be an irreducible Laurent polynomial that is non-degenerate with
respect to its Newton polygon $\Delta(f)$.
Then there exists a canonical divisor $K_{\Delta(f)}$ on (the complete non-singular model of) $U(f)$ for which a basis of the Riemann-Roch space
$\mathcal{L}(K_{\Delta(f)})$ is given by
\[ \left\{ \left. \, x^iy^j \, \right| \, (i,j) \in \Delta^{(1)} \cap \mathbb{Z}^2 \right\}.\]
In this, the function field $\mathbb{C}(U(f))$ is understood to be identified with the fraction field of $\mathbb{C}[x^{\pm 1}, y^{\pm 1}] / (f)$.}\\

In particular, this says that the canonical model of $U(f)$ is contained in $\text{Tor}(\Delta^{(1)})$.
The following observation is due to Koelman.

\begin{lemma} \label{gon2char}
Let $f \in \mathbb{C}[x^{\pm 1}, y^{\pm 1}]$ be non-degenerate with respect to its Newton polygon $\Delta(f)$.
Suppose that $\Delta(f)$ is of genus $g \geq 2$. Then $U(f)$ is hyperelliptic if and only if $\Delta(f)^{(1)}$
is one-dimensional.
\end{lemma}

\noindent \textsc{Proof.} See \cite[Lemma~3.2.9]{Koelman}. An alternative proof was
given in \cite[Lemma~5.1]{CastryckVoight} and uses the above reformulation of
Theorem~\ref{Khovanskii}: the function field of the canonical image is $\mathbb{C}(x,y)$ if and only
if $\Delta^{(1)}$ is two-dimensional. \hfill $\blacksquare$\\

As a corollary, we obtain a proof of Conjecture~\ref{gongeneric} in case $\text{lw}(\Delta^{(1)}) = 1$.
Indeed, the above ensures that $U(f)$ generically defines a curve of gonality at least $3$.
But by Theorem~\ref{lwDelta1}, either $\text{lw}(\Delta) = 3$ or $\Delta \cong 4\Sigma$,
and the statement follows from Theorem~\ref{goninequality}.
Entirely similarly, we obtain a proof for the case $\Delta \cong 2 \Upsilon$.
Again, all of this can be turned into an `if and only if'.

\begin{lemma} \label{gon3char}
Let $f \in \mathbb{C}[x^{\pm 1}, y^{\pm 1}]$ be non-degenerate with respect to its Newton polygon
$\Delta(f)$, which we assume to be two-dimensional.
Then $U(f)$ is trigonal if and only if
\[ \emph{lw}(\Delta(f)^{(1)}) = 1 \quad \text{or} \quad \Delta(f) \cong 2 \Upsilon.\]
\end{lemma}

\noindent \textsc{Proof.} It remains to prove the `only if' part.
Parts of the following reasoning already
appeared in an unpublished addendum to \cite{CastryckVoight}. Suppose that
$U(f)$ is trigonal. Then by Petri's theorem, the intersection of all quadrics
containing the canonical image is a rational normal scroll $S \subset \mathbb{P}^{g-1}$, which is a surface
of sectional genus $0$.
On the other hand, by Theorem~\ref{Khovanskii}, $U(f)$
is canonically embedded in $\Tor(\Delta(f)^{(1)})$. An earlier mentioned result by Koelman \cite{Koelman2} states
that $\Tor(\Delta(f)^{(1)})$ is generated by quadrics as soon as $\Delta(f)^{(1)}$ contains
at least $4$ lattice points on the boundary. So:
\begin{itemize}
  \item If $\partial \Delta(f)^{(1)} \geq 4$, then we must have that $\Tor(\Delta(f)^{(1)}) = S$. Since
  $S$ is a surface of sectional genus zero, $\Delta(f)^{(1)}$ cannot have any interior lattice points.
  In particular, either $\text{lw}(\Delta(f)^{(1)}) = 1$, or $\Delta(f)^{(1)} \cong 2 \Sigma$.
  The latter is impossible, however, since then $U(f)$ would be isomorphic to a plane quintic, which has gonality $4$ by
  a result of Namba \cite{Namba} -- see also Theorem~\ref{Namba} below.
  \item If $\partial \Delta(f)^{(1)} = 3$ and $\Delta(f)^{(1)}$ contains an
  interior lattice point, then using Lemma~\ref{haaseschicho}, it is an easy exercise
  to show that $\Delta(f)^{(1)} \cong \Upsilon$, hence $\Delta(f) \cong 2\Upsilon$.
\end{itemize} This concludes the proof. \hfill $\blacksquare$\\

Note that the above proof gives a prudent indication of the exceptionality of $2\Upsilon$.

Again, entirely similarly to the foregoing cases, we can use this to obtain a proof of
Conjecture~\ref{gongeneric} in case $\text{lw}(\Delta^{(1)}) = 2$. We conclude:

\begin{theorem} \label{smallwidth}
Conjecture~\ref{gongeneric} is true for all lattice polygons $\Delta$ for which $\emph{lw}(\Delta^{(1)}) \leq 2$.
\end{theorem}

To continue this type of iteration, we would need if-and-only-if statements for $\text{lw}(\Delta^{(1)}) = 2,3,4, \dots$, similar
to Lemmata~\ref{gon2char} and \ref{gon3char}. In pursuing this, one naturally bumps into Green's canonical conjecture \cite{Green}, which is an
unproven
generalization of Petri's theorem. It states that the Clifford index of $U(f)$ is the smallest integer $p$ for which the canonical ideal of $U(f)$
does not satisfy
\emph{property $N_p$}. The latter is a certain non-vanishing property of the Betti numbers appearing in a minimal free resolution of the ideal -- see \cite[Chapter~9]{Eisenbud} for an introduction.
We were not able to unveil a connection between property $N_p$ and the combinatorics of the Newton polygon.

However, in an attempt to discover such a connection, we have carried out the following experiment, which
provides evidence for Conjecture~\ref{gongeneric} up to genus $13$. This being
ongoing research, we will be concise here.
Up to equivalence, we have enumerated all two-dimensional lattice polygons $\Delta^{(1)}$ that are interior to
a bigger lattice polygon $\Delta$ and that contain between $3$ and $13$ lattice points. There are $176$ such polygons. For
each of these, we have picked a `generic' Laurent polynomial with Newton polygon $\Delta^{(1)(-1)}$.
For each such polynomial, we have computed the Betti table of the corresponding canonical ideal. We have worked over the finite
field $\mathbb{F}_{10007}$ to speed up the computation; this is not expected to influence the outcome.
Our most notable observations
thus far are:
\begin{itemize}
  \item in each of these $176$ cases, Green's canonical conjecture was consistent with Conjecture~\ref{gongeneric};
  \item the following table is missing in Schreyer's conjectured list of Betti tables appearing in genus $10$ \cite[Section~6]{Schreyertables}:
  \[ \begin{array}{ccccccccc} 1 & . & . & . & . & . & . & . & . \\ . & 28 & 105 & 168 & 154 & 70 & 6 & . & . \\
     . & . & 6 & 70 & 154 & 168 & 105 & 28 & . \\ . & . & . & . & . & . & . & . & 1 \\ \end{array} \]
  (it appeared for the six $\Delta^{(1)}$ for which $\# (\Delta^{(1)(1)} \cap \mathbb{Z}^2) = 2$).
\end{itemize}
All computations were carried out using \textsc{Magma} \cite{Magma}.

\section{Reinterpretation of some previous results} \label{previousresults}

In this section we give additional support for Conjecture~\ref{gongeneric} by reinterpreting some previously
obtained results.

\begin{theorem}[Namba, 1979] \label{Namba}
Conjecture~\ref{gongeneric} holds if $\Delta \cong d \Sigma$ for some integer $d \geq 0$.
\end{theorem}

\noindent \textsc{Proof.} If $\Delta \cong d \Sigma$, then a
generic $\Delta$-supported Laurent polynomial $f \in \mathbb{C}[x^{\pm 1}, y^{\pm 1}]$
defines a smooth curve of degree $d$. A result of Namba \cite{Namba} states that such curves have gonality
$d-1$. \hfill $\blacksquare$\\

The above can be generalized to the case where $\Delta^{(1)} \cong d \Sigma$ for some
integer $d \geq 0$ (corresponding to smooth plane curves with, possibly, some prescribed behavior at
the coordinate points).

\begin{theorem}[Martens, 1996] \label{Hirzebruch}
Let $a,b \geq 1$ and $k \geq 0$ be integers.
Then Conjecture~\ref{gongeneric} holds if
\[ \Delta \cong \emph{Conv} \{ (0,0), (a+bk,0),(a,b), (0,b) \}. \]
\end{theorem}

\noindent \textsc{Proof.} Note that $\text{Tor}(\Delta) \cong H_k$, the Hirzebruch surface
of invariant $k$.
If $f \in \mathbb{C}[x^{\pm 1}, y^{\pm 1}]$ is non-degenerate
with respect to its Newton polygon $\Delta$, then $U(f)$ will embed smoothly in $H_k$. Martens \cite{Martens} proved
that the gonality is then computed by a ruling of $H_k$. The ruling is unique if $k \geq 1$ and is given by
vertical projection, which is of degree $b = \text{lw}(\Delta)$. If $k=0$, then there are two
rulings, namely horizontal projection and vertical projection. These are of degree $a$ and $b$, respectively,
and since $\text{lw}(\Delta) = \min \{a,b\}$, the result follows. \hfill $\blacksquare$\\

Note that the case $k=0$, corresponding to rectangular polygons, already follows
from older work of Schreyer \cite{Schreyer}, who studied the gonality of curves in $\text{Tor}(\Delta) \cong \mathbb{P}^1 \times \mathbb{P}^1$.

Again, Theorem~\ref{Hirzebruch} can be adapted to the case where actually
\[ \Delta^{(1)} \cong \text{Conv} \{ (0,0), (a+bk,0),(a,b), (0,b) \}. \]

Recently, Martens' result on Hirzebruch surfaces was generalized to certain
of their blow-ups. This gives the following result, which subsumes Theorem~\ref{Hirzebruch}.

\begin{theorem}[Kawaguchi, 2008, 2010] \label{Kawaguchi}
Let $a,b\geq 1$ and $k \geq 0$ be integers. Let $C \subset \mathbb{R}^2$ be the graph of
a concave, continuous, piece-wise linear function $f : [0,a + bk] \rightarrow \mathbb{R}^+$ with $f(0) > 0, f(a) = b,$ and $kf(a + bk) = 0$,
such that its segments have lattice points as endpoints, there is at least
one horizontal segment, and $f$ is linear on $[a, a+bk]$. Then Conjecture~\ref{gongeneric}
holds for the convex hull of $C \cup \{(0,0),(a,0)\}$.\\
\begin{figure}[ht]
\centering
\includegraphics[height=2cm]{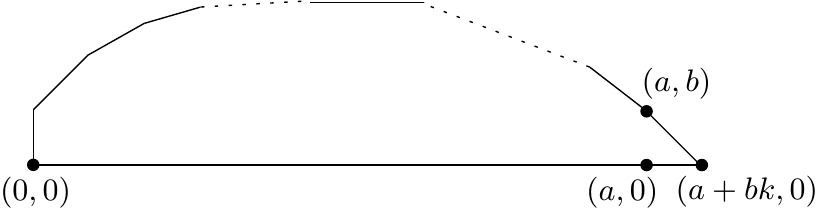}
\end{figure}
\end{theorem}

\noindent \textsc{Proof.} Same proof, now using \cite{Kawaguchi,Kawaguchi2} instead of \cite{Martens}. \hfill $\blacksquare$.\\

We end with a very particular case. To us, it is nevertheless important
because it proves Conjecture~\ref{gongeneric} in a situation where the combinatorial properties are somewhat harder to grasp.

\begin{theorem}[Martens, 1980] \label{Martens}
Conjecture~\ref{gongeneric} holds if $\Delta \cong 3\Upsilon$.
\end{theorem}

\noindent \textsc{Proof.}
$\text{Tor}(3\Upsilon) \cong \text{Tor}(\Upsilon)$ can be realized in $\mathbb{P}^3 = \text{Proj} \, \mathbb{C}[x,y,z,w]$
as the zero locus of $xyz - w^3$. An irreducible Laurent polynomial $f$ that is supported
on $3 \Upsilon$ then defines another cubic in $\mathbb{C}[x,y,z,w]$, the intersection of whose
zero locus with $\text{Tor}(\Upsilon)$ is birationally equivalent to $U(f)$. Generically, this
intersection will be smooth. It is well-known that smooth complete intersections of two cubics in $\mathbb{P}^3$
have gonality $6$. See \cite{ELMS,Martens2}. \hfill $\blacksquare$

\section{Toric degenerations} \label{toricgeometry}

The material in this section is inspired by \cite[Section~2.3]{Shustin}.
It partly extends Section~\ref{toricsurfaces}, but we will be slightly more concise here.

We first recall, in arbitrary dimension $n \geq 1$, the construction of a \emph{toric variety} $\Tor(\Delta)$ associated to a
\emph{lattice polytope} $\Delta \subset \mathbb{R}^n$, i.e.\ the convex hull of a finite number of points of $\mathbb{Z}^n$.
We will assume throughout that $\Delta$ is \emph{very ample}.
The latter is a technical notion for which we refer to \cite{BrunsGubeladze} (it guarantees that $\text{Tor}(\Delta)$, as constructed below, is isomorphic to the abstract toric variety associated to the normal fan $\Sigma(\Delta)$ of $\Delta$, see below), but we note that very-ampleness is
implied by the existence of a subdivision into unimodular simplices (simplices of volume $1/r!$, with $r$ the dimension of $\Delta$). This is automatic if $n\leq2$.
Define $S := \Delta \cap \mathbb{Z}^n$ and consider the map $\phi$ sending a point $(x_1,\ldots,x_n)$ in the $n$-dimensional torus $\mathbb{T}^n$ over $\mathbb{C}$ to the point in $\mathbb{P}^{|S|-1}$ with projective coordinates $(x_1^{i_1}\cdot\ldots\cdot x_n^{i_n})_{(i_1,\ldots,i_n)\in S}$. The Zariski closure of the image of $\phi$ is the toric variety $\Tor(\Delta)$. If $X_{i_1,\ldots,i_n}$ denotes the projective coordinate of $\mathbb{P}^{|S|-1}$ corresponding to $(i_1,\ldots,i_n)\in S$, then all binomials
of the form
$$ \prod_{k=1}^s X_{i_1^{(k)},\ldots,i_n^{(k)}} - \prod_{k=1}^s X_{j_1^{(k)},\ldots,j_n^{(k)}}$$
for which $$ \sum_{k=1}^s (i_1^{(k)},\ldots,i_n^{(k)}) =
\sum_{k=1}^s (j_1^{(k)},\ldots,j_n^{(k)})$$ are zero on $\Tor(\Delta)$. These binomials generate the homogeneous ideal of $\Tor(\Delta)$.
As in Section~\ref{toricsurfaces}, the faces $\tau \subset \Delta$ naturally decompose $\Tor(\Delta)$
in a disjoint union of toric orbits $O(\tau) \cong \mathbb{T}^{\dim \tau}$. One has $\phi(\mathbb{T}^n) = O(\Delta)$,
and one can show that $\text{Tor}(\Delta)$ is a normal variety.

Now let $\Delta \in \mathbb{R}^2$ be a two-dimensional lattice polygon.
Let $\Delta_1, \dots, \Delta_r$ be a \emph{regular subdivision} of $\Delta$, i.e.\ a
collection of two-dimensional lattice polygons for which
there is an upper-convex piece-wise linear function $v: \Delta \to \mathbb{R}$ such that
  $\Delta_1, \dots, \Delta_r$ are the maximal closed subsets on which $v$ is linear.
If such a $v$ exists, we may assume that $v(i,j)\in\mathbb{Z}$ for each $(i,j)\in\Delta\cap\mathbb{Z}^2$ --
see \cite[Proposition~1.69(i)]{BrunsGubeladze}. Let $\delta$ be an integer
such that $\delta$ is strictly greater than each of these $v(i,j)$'s and
let $\tilde{\Delta}$ be the convex hull in $\mathbb{R}^3$ of all the points
$(i,j,v(i,j))$ and $(i,j,\delta)$ with $(i,j) \in \Delta \cap \mathbb{Z}^2$.
The latter is easily seen to be very ample.
For $\ell \in \{1, \dots, r\}$, we let $\tilde{\Delta}_\ell$ be the face
\[ \{\, (i,j,v(i,j)) \, | \, (i,j) \in \Delta_\ell \, \} \subset \tilde{\Delta}.\]

Let $\tilde{S}=\tilde{\Delta}\cap\mathbb{Z}^3$ and consider the toric threefold $Y=\Tor(\tilde{\Delta})$ in $\mathbb{P}^{|\tilde{S}|-1}$, along
with the corresponding monomial map $\tilde{\phi}:\mathbb{T}^3\hookrightarrow\mathbb{P}^{|\tilde{S}|-1}$.
There is a natural fibration $$p:Y\to \mathbb{P}^1:P=(X_{i,j,k})_{(i,j,k)\in \tilde{S}}\mapsto p(P)=(X_{i,j,k+1}:X_{i,j,k})$$
where $i,j,k$ are chosen such that $X_{i,j,k}$ and $X_{i,j,k+1}$ are not both zero.
The image $p(P)$ is independent of this choice, and for $(x,y,t)\in\mathbb{T}^3$ one has $p(\tilde{\phi}(x,y,t))=(t:1)$.
The fiber $Y_{\infty}:=p^{-1}(1:0)$ is equal to the copy of $\Tor(\Delta)$ contained in the linear subspace $V$ of $\mathbb{P}^{|\tilde{S}|-1}$, defined by $X_{i,j,k}=0$ for all $(i,j,k)\in \tilde{S}$ with $k<\delta$ (i.e., $Y_{\infty}$ is the toric orbit
associated to the top face of $\tilde{\Delta}$). If $t\in \mathbb{C}\setminus\{0\}$, the restriction of the projection $$\pi:\mathbb{P}^{|\tilde{S}|-1}\to V:(X_{i,j,k})_{(i,j,k)\in \tilde{S}}\mapsto (X_{i,j,\delta})_{(i,j)\in\Delta\cap\mathbb{Z}^2}$$ to $Y_t:=p^{-1}(t:1)$ is an isomorphism between $Y_t$ and $Y_{\infty}$.
On the other hand, the fiber $Y_0:=p^{-1}(0:1)$ is equal to $\bigcup_{\ell=1}^r \Tor(\tilde{\Delta}_{\ell})\cong\bigcup_{\ell=1}^r \Tor(\Delta_{\ell})$. So we get a degeneration of the toric surface $\Tor(\Delta)$ to $\bigcup_{\ell=1}^r \Tor(\Delta_{\ell})$.

Let $R = \mathbb{C}[t]$ be equipped with the natural $t$-adic valuation $\text{val} : R \setminus \{0\} \rightarrow \mathbb{Z}$. Let
$$f_t = \sum_{(i,j)\in\Delta\cap\mathbb{Z}^2} a_{i,j}(t)x^iy^j$$
be a Laurent polynomial with coefficients in $R$ that is supported on $\Delta$, such that, when considered
as a trivariate polynomial over $\mathbb{C}$, it is supported on $\tilde{\Delta}$. In particular,
for all $(i,j) \in \Delta \cap \mathbb{Z}^2$ one has $\text{val}\, a_{i,j}(t) \geq v(i,j)$. Define
$c_{i,j} = (a_{i,j}(t) \cdot t^{-v(i,j)}) |_{t=0}$. We make two assumptions about $f_t$:
\begin{itemize}
  \item $f_t$ is non-degenerate with respect to $\Delta$ (when considered as a Laurent polynomial over the
  field of Puiseux series $\mathbb{C}\{\{t\}\}$): this will be referred to as the \emph{non-degeneracy} of $f_t$;
  \item the Laurent polynomials
  \[ \sum_{(i,j) \in \Delta_\ell \cap \mathbb{Z}^2} c_{i,j} x^iy^j \ \ \in \ \mathbb{C}[x^{\pm 1}, y^{\pm 1}] \]
  (for $\ell = 1, \dots, r$) are non-degenerate with respect to the respective $\Delta_\ell$: this will be referred to
  as the \emph{local non-degeneracy} of $f_t$.
\end{itemize}
Now, when considered as an element of $\mathbb{C}[t,x,y]$, our Laurent polynomial $f_t$ defines a hyperplane section $X$ of $Y$. For $t\in\mathbb{C}\setminus\{0\}$, the fiber $X_t:=X\cap Y_t$ is equal to the intersection of $Y_t$ with the hyperplane defined by
$$\sum_{(i,j)\in\Delta\cap\mathbb{Z}^2} a_{i,j}(t)X_{i,j,\delta}=0.$$
The fiber $X_0:=X\cap Y_0$ of $X$
is equal to $$X\cap \bigcup_{\ell=1}^r \Tor(\tilde{\Delta}_{\ell})=\bigcup_{\ell=1}^r (X\cap \Tor(\tilde{\Delta}_{\ell})),$$ where $X^{(\ell)}:=X\cap \Tor(\tilde{\Delta}_{\ell})$ is the intersection of $\Tor(\tilde{\Delta}_{\ell})$ with the hyperplane
$$\sum_{(i,j)\in\Delta_{\ell}\cap\mathbb{Z}^2} c_{i,j} X_{i,j,v(i,j)}=0.$$
A pair $X^{(\ell)},X^{(m)}$ intersects if and only if the polygons $\tilde{\Delta}_\ell$ and $\tilde{\Delta}_m$ have an edge in common, and if so, the intersection is defined by the hyperplane section $$\sum_{(i,j,k)\in\tilde{\Delta}_{\ell}\cap\tilde{\Delta}_m\cap\mathbb{Z}^3} c_{i,j} X_{i,j,k}=0$$ of $\Tor(\tilde{\Delta}_{\ell})\cap\Tor(\tilde{\Delta}_m)=\Tor(\tilde{\Delta}_{\ell}\cap\tilde{\Delta}_m)$.
By local non-degeneracy, the curves $X^{(\ell)}$ are smooth. Moreover, $X^{(\ell)}$ and $X^{(m)}$
intersect $\Tor(\tilde{\Delta}_{\ell}\cap\tilde{\Delta}_m)$ transversally in the same points. Hence
they intersect each other transversally. The number of intersection points
is equal to the number of lattice points in $\tilde{\Delta}_{\ell}\cap\tilde{\Delta}_m$ minus one. For instance, if $\tilde{\Delta}_{\ell}\cap\tilde{\Delta}_m$ is a line segment without lattice points in its interior, then the curve $\Tor(\tilde{\Delta}_{\ell}\cap\tilde{\Delta}_m)$ is a projective line and $X^{(\ell)}\cap X^{(m)}$ is a point on this line.

What we have actually constructed is a \emph{strongly semi-stable arithmetic surface} over $\mathbb{C}[[t]]$ (see \cite[Section~1.1]{MBaker}
for this terminology). Indeed, consider the restriction of $p$ to $X^\text{fin} := X \setminus X_\infty ( = X \setminus p^{-1}(1:0))$. This gives $X^\text{fin}$ the structure of a scheme
over $\mathbb{A}^1 = \text{Spec} \, \mathbb{C}[t]$. Define $\mathfrak{X} := X^\text{fin} \otimes \mathbb{C}[[t]]$
as a scheme over $\mathbb{C}[[t]]$. It is proper and flat, and its generic fiber $\mathfrak{X} \otimes \mathbb{C}\{\{t\}\}$
is precisely the Zariski-closed embedding of
$U(f_t)$ in $\text{Tor}(\Delta)$, as described in Section~\ref{toricsurfaces} (with $\mathbb{C}$ replaced by $\mathbb{C}\{\{t\}\}$). By non-degeneracy, this is a smooth curve,
hence $\mathfrak{X}$ is an arithmetic surface.
On the other hand, the special fiber $\mathfrak{X} \otimes \mathbb{C}$ is precisely the reducible curve
$X_0$ having $X^{(1)}, \dots, X^{(r)}$ as its components. By local non-degeneracy, these components are smooth and intersect each other
transversally. Hence the reduction is strongly semi-stable.\\

\noindent \textbf{Example.} Let $f_t=1+x+y+txy$ and let $\Delta$, $\Delta_1$, $\Delta_2$ and $\tilde{\Delta}$ be as depicted below:

\begin{figure}[h]
\centering
\subfloat{\includegraphics[height=3cm]{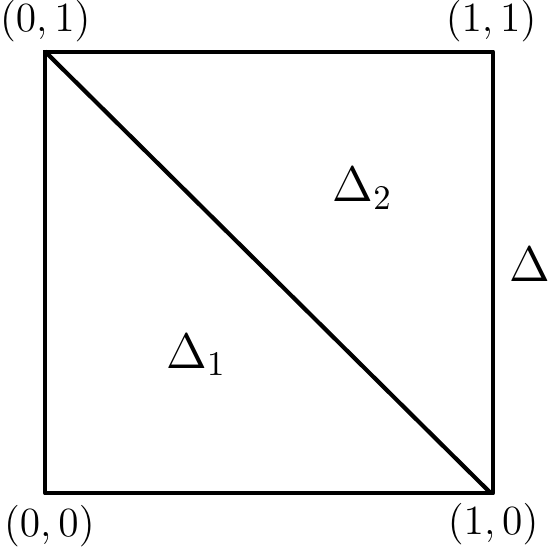}}
\hspace{1cm}
\subfloat{\includegraphics[height=3cm]{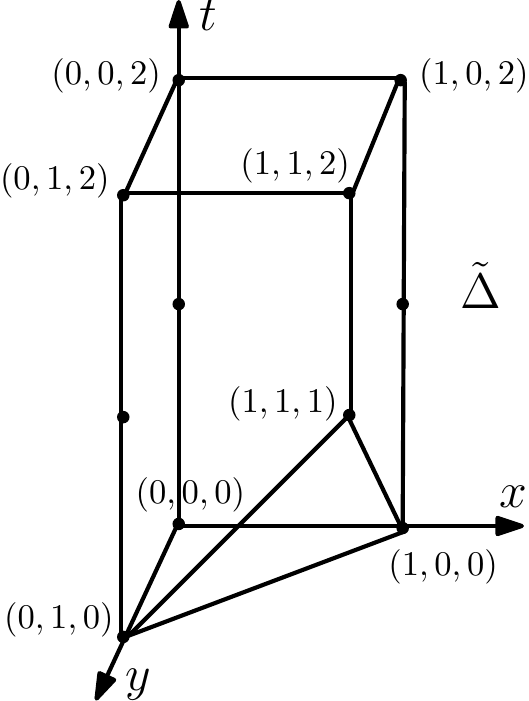}}
\end{figure}
\noindent and $$v:(x,y)\in\Delta\mapsto \begin{cases} 0 & \text{if}\ (x,y)\in\Delta_1 \\ x+y-1 & \text{if}\ (x,y)\in\Delta_2. \end{cases}$$ The toric threefold $Y=\Tor(\tilde{\Delta})$ lies in $\mathbb{P}^{10}$. For $t \in \mathbb{C} \setminus \{0\}$, the fiber $Y_t$ is isomorphic to $\Tor(\Delta)=\mathbb{P}^1\times\mathbb{P}^1$ and its special fiber $Y_0$ is the union of two planes $\Tor(\tilde{\Delta}_1)$ and $\Tor(\tilde{\Delta}_2)$ (respectively with projective coordinates $(X_{0,0,0}:X_{1,0,0}:X_{0,1,0})$ and $(X_{1,0,0}:X_{0,1,0}:X_{1,1,1})$) that intersect each other in the line $\Tor(\tilde{\Delta}_1\cap\tilde{\Delta}_2)$ (with projective coordinates $(X_{1,0,0}:X_{0,1,0})$). The special fiber $X_0$ is the union of two lines, namely $X^{(1)}\subset \Tor(\tilde{\Delta}_1)$ with equation $X_{0,0,0}+X_{1,0,0}+X_{0,1,0}=0$ and $X^{(2)}\subset \Tor(\tilde{\Delta}_2)$ with equation $X_{1,0,0}+X_{0,1,0}+X_{1,1,1}=0$, that intersect in the point $(1:-1)\in\Tor(\tilde{\Delta}_1\cap\tilde{\Delta}_2)$.\\

For the application that we have in mind, our strongly semi-stable arithmetic surface $\mathfrak{X}$ is supposed to be \emph{regular}, which
in our case is equivalent to saying that the singular locus of $X$ does not meet $X_0$. In general, this is not satisfied. However, by local non-degeneracy,
the singularities at $X_0$ are entirely related to the fact that the ambient space $Y = \text{Tor}(\tilde{\Delta})$ is itself
singular at $Y_0$, and a \emph{toric resolution} automatically resolves the singularities of $X$ at $X_0$.
We give a brief sketch, in which we assume some additional background concerning toric
varieties $\text{Tor}(\Sigma)$ constructed from fans $\Sigma$. For an account on
such abstract toric varieties and toric resolutions, see \cite{cox}. For more
details on resolving non-degenerate hypersurface singularities, see \cite{kouchnirenko}.

Let $\Sigma(\tilde{\Delta})$ be the normal fan of $\tilde{\Delta}$. One can always find a subdivision $\Sigma'$ of $\Sigma(\tilde{\Delta})$ such that
the induced birational morphism $\rho: \text{Tor}(\Sigma') \rightarrow \text{Tor}(\Sigma(\tilde{\Delta})) \cong \text{Tor}(\tilde{\Delta})$ is
a resolution of singularities.
Write $Y' = \text{Tor}(\Sigma')$ and let
$X' \subset Y'$ be the strict transform of $X$ under $\rho$. The morphism
$p' = p \circ \rho$ yields a fibration $Y' \rightarrow \mathbb{P}^1$. One can then redo the argument and obtain an arithmetic surface $\mathfrak{X}'$ over $\mathbb{C}[[t]]$,
which is still strongly semi-stable, but which is moreover regular.
The generic fibers of $\mathfrak{X}$ and $\mathfrak{X}'$ are isomorphic, because $\rho|_{X'}$ is an isomorphism on $p'^{-1}(V)$
for an open subset $V$ of $\mathbb{P}^1$.
On the other hand, the special fiber of $\mathfrak{X}'$ differs from the special fiber of $\mathfrak{X}$.
To see how the latter modifies under toric resolutions, it suffices to analyze what
happens when we subdivide a two-dimensional cone.

First, we consider cones $\sigma^2_{\ell,m}$ spanned
by rays $\sigma^1_{\ell}$ and $\sigma^1_m$ that correspond to adjacent lower facets
$\tilde{\Delta}_\ell$ and $\tilde{\Delta}_m$ of $\tilde{\Delta}$.
Then $\sigma^2_{\ell,m}$ corresponds to the edge $\tilde{\Delta}_\ell \cap \tilde{\Delta}_m$, and
the introduction of a new
ray boils down to blowing up $Y$ in $\text{Tor}(\tilde{\Delta}_\ell \cap \tilde{\Delta}_m)$. This separates the curves
$X^{(\ell)}$ and $X^{(m)}$, and each intersection point becomes replaced by an exceptional curve intersecting $X^{(\ell)}$ and $X^{(m)}$
transversally. This exceptional curve
is contained in the strict transform of $X$ and hence belongs to the special fiber of our new arithmetic surface. All intersections remain transversal.
\begin{figure}[h]
\centering
\subfloat{\includegraphics[height=3cm]{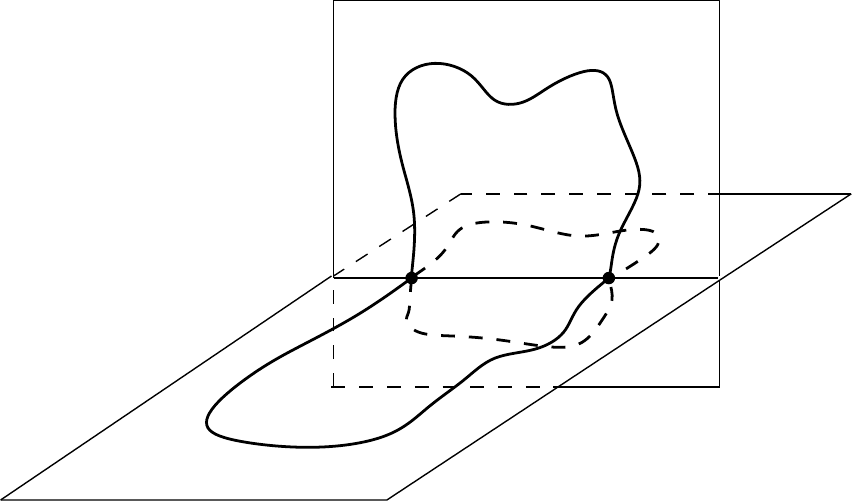}}
\hspace{1cm}
\subfloat{\includegraphics[height=3cm]{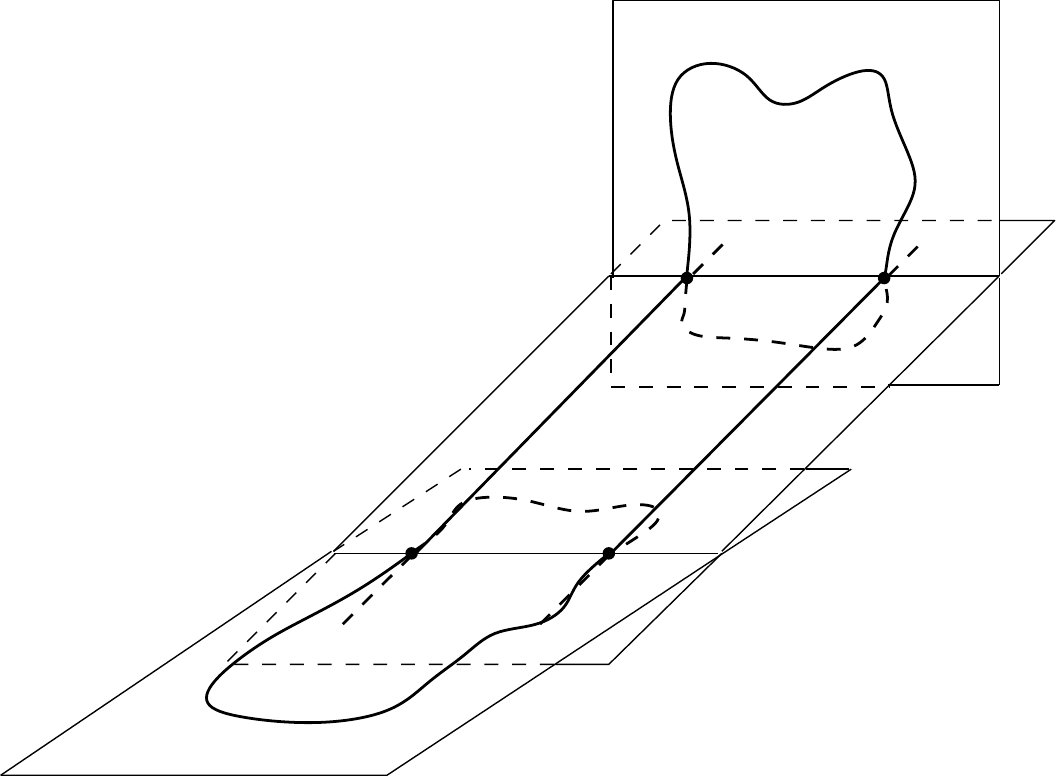}}
\end{figure}
More generally, if $k$ rays are added to $\sigma_{\ell,m}$, then each intersection point becomes replaced by
a chain of $k$ transversally intersecting exceptional curves.

Next, consider cones $\sigma^2_\ell$ spanned by a ray $\sigma^1_{\ell}$ that corresponds to a lower facet $\tilde{\Delta}_\ell$
and a ray that corresponds to an adjacent vertical facet $\tilde{\Gamma}$. Then the introduction of a new ray boils down to blowing up $Y$ in
$\text{Tor}(\tilde{\Delta}_\ell \cap \tilde{\Gamma})$. Each point of intersection of $X^{(\ell)}$ with $\text{Tor}(\tilde{\Delta}_\ell \cap \tilde{\Gamma})$ becomes equipped with an emanating exceptional curve. More generally, if $k$ rays are added to $\sigma^2_\ell$,
then the point becomes equipped with an emanating chain of $k$ transversally intersecting exceptional curves.

\section{The gonality of the dual graph} \label{graphgonalities}

Linear systems on graphs and on metric graphs have been introduced by Baker and Norine in \cite{MBaker,BakNor}.
It turns out that these linear systems obey properties that are analogous to those of linear systems on algebraic curves, such as the Riemann-Roch Theorem. Moreover, the Specialization Lemma \cite[Lemma 2.8]{MBaker} can be used to transport results on metric graphs to algebraic curves or vice versa. In this section, we will
briefly overview the basic notions of linear systems on metric graphs and use the Specialization Lemma to obtain a lower bound for curve gonalities.

Let $G=(V(G),E(G))$ be a connected graph without loops. The {\it metric graph} $\Gamma$ assiociated to $G$ is the compact connected metric space where each edge is identified with a line segment $[0,1]$. A {\it divisor} $D=a_1v_1+\ldots+a_sv_s$ on $\Gamma$ is an element of the free abelian group generated by the points of the metric graph $\Gamma$. The {\it degree} of $D$ is the sum $a_1+\ldots+a_s$ of its coefficients and $D$ is {\it effective} if and only if each coefficient $a_i$ is nonnegative. Let $\phi:\Gamma\to\mathbb{R}$ be a continuous map such that the restriction of $\phi$ to an edge of $\Gamma$ is piece-wise linear with integer slopes and only finitely many pieces. If $v\in\Gamma$, write $\text{ord}_v(\phi)$ to denote the sum of the incoming slopes of $\phi$ at $v$. Note that $\text{ord}_v(\phi)$ is nonzero for only finitely many points $v\in\Gamma$, so we can consider the divisor $\text{div}(\phi)=\sum_{v\in\Gamma} \left(\text{ord}_v{\phi}\right) \cdot v$. We say that two divisors $D$ and $D'$ on $\Gamma$ are {\it equivalent}, and denote this by $D\sim D'$, if and only if $D'-D=\text{div}(\phi)$ for some $\phi$. The {\it complete linear system} $|D|$ of a divisor $D$ is the set of all effective divisors $D'$ that are equivalent to $D$. The {\it rank} $r(D)$ of the linear system $|D|$ is defined as follows. We have that $r(D)=-1$ if and only if $|D|=\emptyset$ and $r(D)\geq r$ if and only if $|D-E|\neq\emptyset$ for all effective divisors $E$ on $\Gamma$ of degree $r$. For instance, $r(D)=1$ if and only if $|D-P|\neq \emptyset$ for each point $P\in\Gamma$ and there exist points $P_1,P_2\in\Gamma$ (not necessarily distinct) such that $|D-P_1-P_2|=\emptyset$. The {\it gonality} of $\Gamma$ is the minimal degree of a divisor on $\Gamma$ having rank one.

For a lattice polygon $\Delta \subset \mathbb{R}^2$ and a regular subdivision $\Delta_1,\ldots,\Delta_r\subset\Delta$, let $G=G(\Delta_1,\ldots,\Delta_r)$ be the graph with vertex set $V(G)=\{v_1,\ldots,v_r\}$ such that the number of edges between the vertices $v_{\ell}$ and $v_m$ is equal to the number of lattice points of $\Delta_{\ell}\cap\Delta_m$ minus one. Let $\Gamma=\Gamma(\Delta_1,\ldots,\Delta_r)$ be the metric graph associated to $G$.

\begin{theorem} \label{graphgon}
Let $\Delta\subset \mathbb{R}^2$ be a two-dimensional
lattice polygon and let $\Delta_1,\ldots,\Delta_r$ be a regular subdivision of $\Delta$.
Let $S$ be the
set of irreducible Laurent polynomials $f \in \mathbb{C}[x^{\pm 1}, y^{\pm 1}]$ for which $\Delta(f) = \Delta$
and the gonality of $U(f)$ is at least the gonality of the metric graph $\Gamma(\Delta_1,\ldots,\Delta_r)$.
Then $S$
is Zariski dense in the space of Laurent polynomials $f \in \mathbb{C}[x^{\pm 1}, y^{\pm 1}]$ for which $\Delta(f) \subset \Delta$.
\end{theorem}

\noindent \textsc{Proof.}
Let $R = \mathbb{C}[t]$.
We construct a Laurent polynomial $f_t \in R[x^{\pm 1}, y^{\pm 1}]$ to which we can apply the
machinery of Section~\ref{toricgeometry}. Let $v : \Delta \rightarrow \mathbb{R}$ be
a un upper-convex piece-wise linear function realizing the subdivision $\Delta_1, \dots, \Delta_r$ such
that $v(\Delta \cap \mathbb{Z}^2) \subset \mathbb{Z}$.

First, let
\[ g_t = \sum_{(i,j) \in \Delta \cap \mathbb{Z}^2} c_{i,j} x^iy^j \quad \in \mathbb{C}[x^{\pm 1}, y^{\pm 1}] \]
be such that each $g_{t, \ell} := \sum_{(i,j) \in \Delta_\ell \cap \mathbb{Z}^2} c_{i,j} x^iy^j$ is non-degenerate with
respect to its Newton polygon $\Delta_\ell$. This is possible, because each non-degeneracy condition is generically satisfied,
and it will guarantee the local non-degeneracy of $f_t$ below. Second, we consider the polynomial
\[ h_t = \sum_{(i,j) \in \Delta \cap \mathbb{Z}^2} c_{i,j} t^{v(i,j)} x^iy^j \quad \in \mathbb{C}[t][x^{\pm 1}, y^{\pm 1}].\]
Now since $\mathbb{C}[t]$ is infinite and since non-degeneracy is generically satisfied, there does exist a Laurent polynomial
\[ h'_t = \sum_{(i,j) \in \Delta \cap \mathbb{Z}^2} a_{i,j}(t) x^iy^j \quad \in \mathbb{C}[t][x^{\pm 1}, y^{\pm 1}]\]
that is non-degenerate with respect to its Newton polygon $\Delta$ (when considered as a Laurent polynomial
with coefficients in $\mathbb{C}\{\{t\}\}$).
But then all but finitely many among the Laurent polynomials
\[ h_t + \lambda(h_t' - h_t), \quad \lambda \in \mathbb{C}[t] \]
must be non-degenerate with respect to their Newton polygon $\Delta$: indeed, this spans a line in coefficient space which
is not strictly contained in the degenerate locus. By taking a $\lambda$ with high $t$-adic valuation,
we end up with a Laurent polynomial $f_t$ that is non-degenerate with respect to its Newton polygon $\Delta$, such that the
$t$-adically leading terms of the coefficients are the same as in $h_t$.

Then by letting $\delta$ be an integer that is strictly
bigger than the valuation of each of the coefficients of $f_t$, and by constructing $\tilde{\Delta}$ accordingly, we
can follow Section~\ref{toricgeometry} and
end up with a (possibly non-regular) strongly semi-stable arithmetic surface $\mathfrak{X}$ over $\mathbb{C}[[t]]$.
The dual
graph of the special fiber $\mathfrak{X} \otimes \mathbb{C} = X^{(1)} \cup \ldots \cup X^{(r)}$
is equal to $G(\Delta_1,\ldots,\Delta_r)$. Indeed, each vertex $v_{\ell}$ corresponds to a curve $X^{(\ell)}$ and each edge $e=(v_{\ell},v_m)$ corresponds to an intersection point of $X^{(\ell)}$ and $X^{(m)}$.
Let $\Gamma = \Gamma(\Delta_1,\ldots,\Delta_r)$ be the associated metric graph.
Now let $\mathfrak{X}'$ be a \emph{regular} strongly semi-stable arithmetic surface, obtained from a subdivision of $\Sigma(\tilde{\Delta})$.
By refining the subdivision if necessary, we may assume that each two-dimensional cone of $\Sigma(\tilde{\Delta})$ becomes
subdivided by an equal amount of rays (say $k$). Then the dual graph of $\mathfrak{X}' \otimes \mathbb{C}$ is obtained
from $G(\Delta_1,\ldots,\Delta_r)$ by introducing $k$ new vertices on each edge, and by attaching
to certain vertices an emanating linear graph. Denote it by $G'$ and let $\Gamma'$ be the associated metric graph.
Then the gonalities of $\Gamma$ and $\Gamma'$ are the same. Indeed, removing the emanating linear graphs from $\Gamma'$ clearly does
not affect the gonality, and the remaining graph is a mere rescaling of $\Gamma$ (by a factor $k+1$).

Then \cite[Corollary 3.2]{MBaker} implies that the gonality of $U(f_t)$ over $\mathbb{C}\{\{t\}\}$ is at least the gonality of $\Gamma'$, hence it
is at least the gonality of $\Gamma$.
To be precise, in \cite{MBaker}, the results are
stated using the $\mathbb{Q}$-graph $\Gamma'_{\mathbb{Q}}$ (i.e.\ only the rational points on the edges are considered), but the gonality of a metric graph $\Gamma'$ is equal to the gonality of its corresponding $\mathbb{Q}$-graph $\Gamma'_{\mathbb{Q}}$. Indeed, by \cite[Corrolary 1.5]{MBaker}, a $\mathbb{Q}$-divisor has rank one on $\Gamma'$ if and only if it has rank one on $\Gamma'_{\mathbb{Q}}$, so the gonality of $\Gamma'$ is at least the gonality of $\Gamma'_{\mathbb{Q}}$. On the other hand, using the rational approximation argument from \cite{GathKerber}, it follows that the gonality of $\Gamma'$ is at most the gonality of $\Gamma'_{\mathbb{Q}}$.

Because $\mathbb{C} \cong \mathbb{C} \{\{t \} \}$, there exists a Laurent polynomial $f' \in \mathbb{C}[x^{\pm 1}, y^{\pm 1}]$ such that the gonality of $U(f')$ over $\mathbb{C}$ is equal to the gonality of $U(f_t)$ over $\mathbb{C}\{\{t\}\}$. Thus the gonality of $U(f')$ is bounded
from below by the gonality of $\Gamma$.
To conclude the proof,
one can either analyze the degree of freedom in the construction of $f_t$ above, or apply the semi-continuity lemma below. \hfill $\blacksquare$

\begin{lemma}[semi-continuity] \label{semcon}
Let $\Delta$ be a lattice polygon. Let
$M_\Delta \subset \mathbb{C}[x^{\pm 1}, y^{\pm 1}]$ be the
set of Laurent polynomials that are non-degenerate with respect to their Newton polygon
$\Delta$ (seen as a quasi-affine variety in coefficient space). Then the map $M_\Delta \rightarrow \mathbb{Z}$ sending $f$ to
the gonality of $U(f)$ is lower semi-continuous.
\end{lemma}

\noindent \textsc{Proof.}
Let $g$ be the genus of $\Delta$ and let $\mathcal{M}_g$ be the moduli space of curves of genus $g$. It is well-known that the map
$\mathcal{M}_g \rightarrow \mathbb{Z}$ sending a curve to its gonality is lower semi-continuous -- see e.g. \cite[Prop.~3.4]{LangeMartens}.
By the flatness of the family of curves parameterized by
$M_\Delta$, we are given a unique morphism $M_\Delta \rightarrow \mathcal{M}_g$ sending $f$ to
the isomorphism class of $U(f)$. See \cite[Section~2]{CastryckVoight} for more details. Since $M_\Delta$ is irreducible, the result follows. \hfill $\blacksquare$\\

We expect that Theorem~\ref{graphgon} is sharp, in the following sense:

\begin{conjecture} \label{conjecturelowerbound}
Let $\Delta\subset\mathbb{R}^2$ be a two-dimensional lattice polygon. Then there is a regular subdivision $\Delta_1,\ldots,\Delta_r$ of $\Delta$ such that the
set of irreducible Laurent polynomials $f \in \mathbb{C}[x^{\pm 1}, y^{\pm 1}]$ for which
$\Delta(f) = \Delta$ and the
gonality of $U(f)$ is equal to the gonality of the metric graph $\Gamma(\Delta_1,\ldots,\Delta_r)$,
is Zariski dense in the space of Laurent polynomials $f \in \mathbb{C}[x^{\pm 1}, y^{\pm 1}]$ with $\Delta(f) \subset \Delta$.
\end{conjecture}

The above conjecture is true for $\Delta=2\Upsilon$, since the metric graph $\Gamma$ corresponding to the subdivision $\Delta_1,\ldots,\Delta_{12}$ of $2\Upsilon$ (see the picture below) has gonality equal to $3$. For instance, the divisor $v_1+v_2+v_3$ has rank one.

\begin{figure}[h]
\centering
\subfloat{\includegraphics[height=3cm]{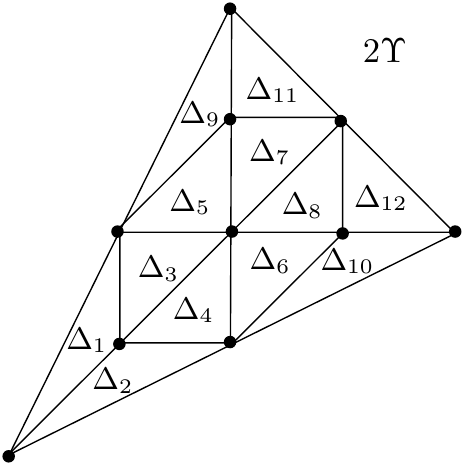}}
\hspace{1,5cm}
\subfloat{\includegraphics[height=3cm]{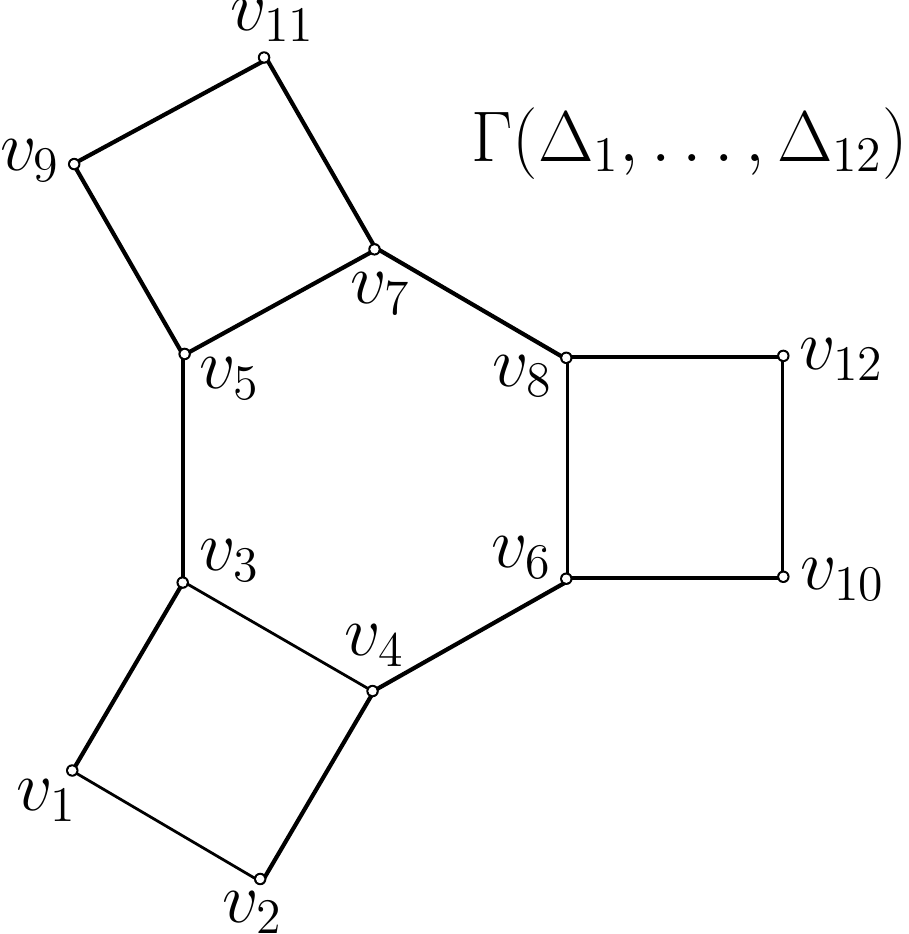}}
\end{figure}

\section{A purely combinatorial conjecture} \label{purelycombinatorial}

Conjecture~\ref{gongeneric} and Conjecture~\ref{conjecturelowerbound} can be combined to a purely combinatorial statement.

\begin{conjecture} \label{conjecturecomb}
Let $\Delta\subset\mathbb{R}^2$ be a two-dimensional lattice polygon. Then there exists a regular subdivision $\Delta_1,\ldots,\Delta_r$ of $\Delta$ such that the gonality of the metric graph $\Gamma(\Delta_1,\ldots,\Delta_r)$ is equal to $\emph{lw}(\Delta^{(1)})+2$ if $\Delta \not \cong 2\Upsilon$, and to $3$ if $\Delta \cong 2\Upsilon$.
\end{conjecture}

We will prove this conjecture for a particular family of lattice polygons (see Theorem \ref{theoremnewfamily}). For this, we need to study the gonality of a certain metric graph. When dealing with linear systems on metric graphs, it is often convenient to view an effective divisor $D=a_1v_1+\ldots+a_sv_s$ as a chip configuration on $\Gamma$ where a stack of $a_i$ chips is placed on the point $v_i$ of $\Gamma$. We will use the chip terminology, the notion of reduced divisors \cite[Theorem 10]{HKN} and Dhar's burning algorithm \cite[Section 2]{Luo} in the following proof.

\begin{lemma} \label{lemmaspecialmetricgraph}
If $r \geq 1$ be an integer and let $G_r$ be the graph defined by $V(G_r)=\{v_1,\ldots,v_r\}$ and $$E(G_r)=\{e_{i,j}=(v_{i-1},v_i)\,|\,i=2,\ldots,r; j=1,\ldots,i\},$$ where the latter
should be seen as a multiset. Then its corresponding metric graph $\Gamma_r$ has gonality equal to $r$.\\

\begin{figure}[h]
\centering
\includegraphics[height=2cm]{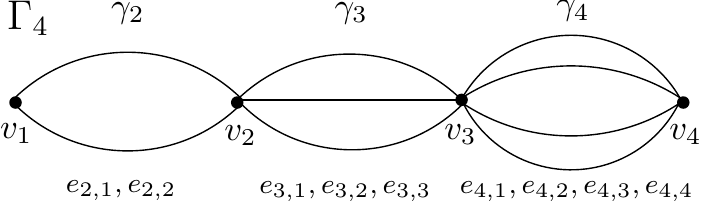}
\end{figure}
\end{lemma}

\noindent \textsc{Proof.}
Since the divisor $v_1+\ldots+v_r$ on $\Gamma$ has rank one, the gonality of $\Gamma$ is at most equal to $r$. Suppose $D$ is an effective divisor on $\Gamma$ with $\deg(D)<r$ and rank at least one. We may assume that $D$ is $v_1$-reduced, hence $D$ has at least one chip at $v_1$. Let $\gamma_i$ be obtained by taking the union of the $i$ edges $e_{i,1},\ldots,e_{i,i}$ between $v_{i-1}$ and $v_i$ and excluding the vertex $v_{i-1}$. Since $\Gamma=\{v_1\}\cup \gamma_2\cup \ldots \cup \gamma_r$, the pigeonhole principle implies that at least one of the subsets $\gamma_2,\ldots,\gamma_r$ of $\Gamma$ does not contain a chip of $D$. Let $i$ be the maximal index for which $\gamma_i$ does not contain a chip of $D$. If we perform Dhar's burning algorithm to reduce $D$ with respect to $v_i$, the chips of $D$ on the subset $\gamma_{i+1}\cup\ldots\cup\gamma_r$ will not move, since $D$ is $v_1$-reduced and hence fire from $v_1$ will pass through $v_i$. So we need that at some point chips must move along $\gamma_i$, since $D$ has rank at least one. If a chip moves along one of the edges of $\gamma_i$, this must also be the case for the other edges. Indeed, otherwise we can find a cycle in $\gamma_i\cup\{v_{i-1}\}\subset \Gamma$ such that chips on it only move in one direction, which cannot happen inside a linear system. We conclude that $D$ must have at least $i$ chips in $\{v_1\}\cup\gamma_2\cup\ldots\cup\gamma_{i-1}$. Since $\gamma_{i+1},\ldots,\gamma_r$ contain at least one chip of $D$, the total amount of chips or the degree of $D$ is at least $i+(r-i)=r$, a contradiction.
\hfill $\blacksquare$

\begin{theorem} \label{theoremnewfamily}
Let $a,b$ be integers with $1\leq a\leq b$. Let $C\subset \mathbb{R}^2$ be the graph of a concave, continuous, piece-wise linear function $f:[0,b]\to\mathbb{R}^+$ with $f(0)=a\geq f(1)$ and $f(b)=0$ such that its segments have lattice points as end points. Then Conjecture \ref{conjecturecomb} holds for the convex hull $\Delta$ of $C$ with $\{(0,0)\}$.
\end{theorem}

\noindent \textsc{Proof.}
Consider the regular subdivision $\Delta_1,\ldots,\Delta_a$ of $\Delta$ where $$\Delta_i=\begin{cases} \text{Conv}\{(0,0),(2,0),(0,2)\} & \text{if } i=1, \\ \text{Conv}\{(i,0),(i+1,0),(0,i+1),(0,i)\}& \text{if } i=2,\ldots,a-1, \\ \Delta\setminus\text{Conv}\{(0,0),(a,0),(0,a)\}& \text{if } i=a.\end{cases}$$

\begin{figure}[h]
\centering
\includegraphics[height=3cm]{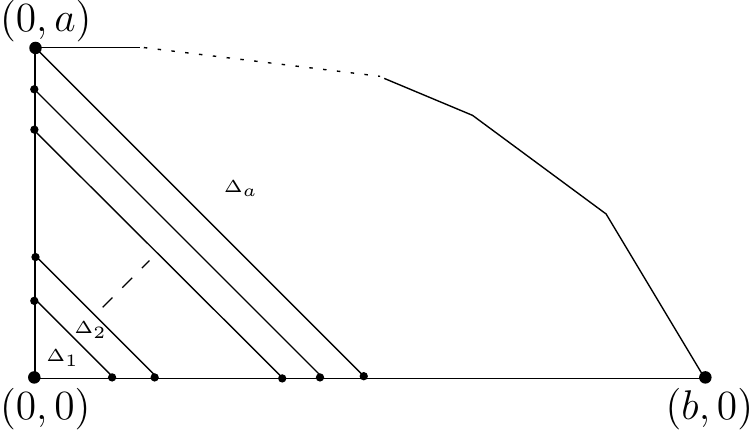}
\end{figure}

If $\Delta \neq a \Sigma$, then $\Gamma(\Delta_1,\ldots,\Delta_a) = \Gamma_a$, which by
Lemma \ref{lemmaspecialmetricgraph} has gonality equal to $a$. On the other hand, $\text{lw}(\Delta^{(1)}) + 2 = \text{lw}(\Delta)$ is equal to $a$
(it cannot be strictly less than $a$ because $\Delta$ has two adjacent edges containing at least $a+1$ lattice points).

If $\Delta = a \Sigma$, then $\Delta_a = \emptyset$ and $\Gamma(\Delta_1,\ldots,\Delta_{a-1}) = \Gamma_{a-1}$, which
has gonality equal to $a-1 = \text{lw}(\Delta^{(1)}) + 2$.
\hfill $\blacksquare$\\

This gives a new proof of Theorem~\ref{Namba} and Theorem~\ref{Hirzebruch}, and parts of Theorem~\ref{Kawaguchi}. But it also
deals with various new cases, including triangles of the form $\text{Conv} \{ (0,0),(b,0),(0,a) \}$.

\noindent \textsc{Katholieke Universiteit Leuven, Departement Wiskunde, Afdeling Algebra, Celestijnenlaan 200, 3001 Leuven (Heverlee), Belgium}\\
\noindent \emph{E-mail addresses:} \verb"wouter.castryck@gmail.com", \verb"filip.cools@wis.kuleuven.be"

\newpage

\begin{center}
  \LARGE{Erratum to `Newton polygons and curve gonalities'}
\end{center}

\normalsize

\subsubsection*{Erratum}
The following statements involving `the metric graph $\Gamma(\Delta_1, \dots, \Delta_r)$',
\begin{itemize}
  \item \cite[Theorem~10]{polgon}
  \item \cite[Conjecture~2]{polgon}
  \item \cite[Conjecture~3]{polgon},
\end{itemize}
are false. The erratum is remedied by replacing $\Gamma(\Delta_1, \dots, \Delta_r)$
by another metric
graph $\Gamma(v)$, which depends on an upper-convex piece-wise linear function $v : \Delta \rightarrow \mathbb{R}$
realizing the given subdivision $\Delta_1, \dots, \Delta_r$ and satisfying $v(\Delta \cap \mathbb{Z}^2) \subset \mathbb{Z}$.
The construction of $\Gamma(v)$ is discussed in Section~\ref{correctgraph} below.

\subsubsection*{Background to the erratum}
We have made a conceptual error in the construction of our regular strongly semi-stable arithmetic surface $\mathfrak{X}$ over
$\mathbb{C}[[t]]$, as explained in \cite[Section~7]{polgon}.
The error lies in the last part, involving toric resolutions of singularities.
Namely, it has been overlooked that the exceptional curves that are introduced during the resolution may appear with non-trivial multiplicities,
turning $\mathfrak{X}$ non-stable.
Whereas our construction suggested that one can keep blowing-up to an arbitrary extent, one should
be much more careful and blow-up just the `right' number of times:
\begin{itemize}
\item all singularities should become resolved (enough blow-ups),
\item no non-trivial multiplicities should appear (not too many blow-ups).
\end{itemize}
Luckily, this `right' number always exists and can be controlled in a purely combinatorial way.\\

\subsubsection*{Acknowledgements} We owe much to Johannes Nicaise and Wim Veys for their patient help.
We would also like to thank F.W.O.-Vlaanderen for its financial support.

\setcounter{section}{0}
\section{The graph $\Gamma(v)$} \label{correctgraph}

Let $\Delta \subset \mathbb{R}^2$ be a two-dimensional lattice polygon.
Let $ \Delta_1, \dots, \Delta_r \subset \Delta$ be a regular subdivision and let
$v : \Delta \rightarrow \mathbb{R}$ be an upper-convex piece-wise linear function realizing this subdivision.
Assume that $v(\Delta \cap \mathbb{Z}^2) \subset \mathbb{Z}$.
Let $G(\Delta_1, \dots, \Delta_r)$ be the graph with vertex set $\{ v_1, \dots, v_r \}$ such that for all $\ell,m$ the number of edges between
$v_\ell$ and $v_m$ is equal to the integral length $L(\ell,m)$ of $\Delta_\ell \cap \Delta_m$ (i.e.\ the number of lattice points minus one).

Let $G(v)$ be obtained from $G(\Delta_1, \dots, \Delta_r)$ by replacing each such edge with a linear graph of length $d(\ell,m)$.
Here, $d(\ell,m)$ is the greatest common divisor of the $(2 \times 2)$-minors of
\[ \begin{pmatrix} a_{\ell 1} & a_{\ell 2} & 1 \\ a_{m 1} & a_{m 2} & 1 \\ \end{pmatrix}, \]
where $(a_{\ell 1}, a_{\ell 2}, 1)$ and $(a_{m 1}, a_{m 2}, 1)$ are primitive normal
vectors to $v(\Delta_\ell)$ and $v(\Delta_m)$, respectively. The third coordinate can be taken $1$
because $v(\Delta \cap \mathbb{Z}^2) \subset \mathbb{Z}$.

Finally, let $\Gamma(v)$ be the metric graph associated to $G(v)$, obtained by identifying each edge with the unit interval.

\subsubsection*{Example}

Consider $\Delta = \text{Conv} \{ (-3,0), (3,0), (0,3) \}$. Let $v : \Delta \rightarrow \mathbb{R}$ be the piece-wise linear function whose
graph is the lower convex hull of $$\{ (-1,1,0), (1,1,0), (0,2,0), (-3,0,1), (3,0,1), (0,3,1) \}.$$
  \begin{figure}[H]
  \centering
  \includegraphics[height=2cm]{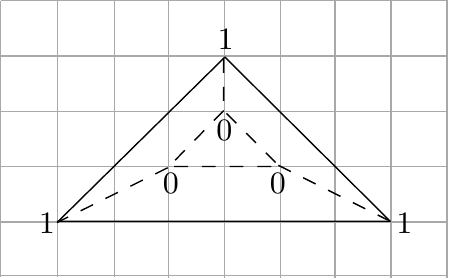}
  \end{figure}
\noindent Denoting the induced subdivision by $\Delta_1, \Delta_2, \Delta_3, \Delta_4$,
one finds that the graph $G(\Delta_1, \Delta_2, \Delta_3 , \Delta_4)$ equals
  \begin{figure}[H]
  \centering
  \includegraphics[height=2cm]{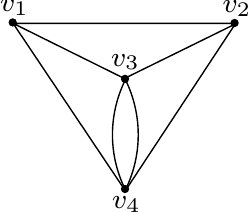}
  \end{figure}
\noindent However, it is easily verified that $G(v)$ equals
  \begin{figure}[H]
  \centering
  \includegraphics[height=2cm]{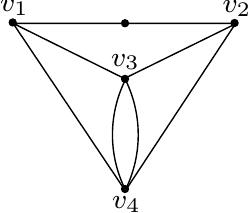}
  \end{figure}

\noindent \emph{Remark.}
By Khovanskii's theorem, a
generic Laurent polynomial $f \in \mathbb{C}[x^{\pm 1}, y^{\pm 1}]$ for which $\Delta(f) \subset \Delta$,
defines a non-hyperelliptic genus $4$ curve $U(f)$ which canonically embeds into $\text{Tor}(\Delta^{(1)}) \subset \mathbb{P}^3$.
Since $\text{Tor}(\Delta^{(1)})$ is a cone, by \cite[Example~5.5.2]{hartshorne} this curve carries a unique $g^1_3$ (computed by
projection from the singular top of the cone).
If $\Gamma(\Delta_1, \Delta_2, \Delta_3, \Delta_4)$ is the metric graph associated to $G(\Delta_1, \dots, \Delta_4)$,
then $\Gamma(\Delta_1, \Delta_2, \Delta_3, \Delta_4)$ carries at least two distinct $g^1_3$'s: it can be verified that
the divisors $3v_1$ and $3v_2$
are non-equivalent. One concludes that $\Gamma(\Delta_1, \Delta_2, \Delta_3, \Delta_4)$ cannot
be the `correct' metric graph associated to this example, since the existence of two distinct $g^1_3$'s would contradict M.\ Baker's specialization theory \cite[Lemma~2.1 and Remark~2.3]{baker}.
In the case of $\Gamma(v)$, the divisors $3v_1$ and $3v_2$ are easily seen to become equivalent.

\section{Details of the toric resolution} \label{proof}

We resume at \cite[Section~7]{polgon}, right after the sentence `For more details on resolving
non-degenerate hypersurface singularities, \dots'\\

Let $\Sigma(\tilde{\Delta})$ be the normal fan of $\tilde{\Delta}$. Any subdivision $\Sigma'$ of
$\Sigma(\tilde{\Delta})$ induces a birational morphism $\rho : \text{Tor}(\Sigma') \rightarrow \text{Tor}(\Sigma(\tilde{\Delta}))
\cong \text{Tor}(\tilde{\Delta})$. If one writes $Y' = \text{Tor}(\Sigma')$ and let $X' \subset Y'$ be the strict transform
of $X$ under $\rho$, then the morphism $p' = p \circ \rho$ yields a fibration $Y' \rightarrow \mathbb{P}^1$, and one can
redo the argument to obtain an arithmetic surface $\mathfrak{X}'$ over $\mathbb{C}[[t]]$.
One can always choose $\Sigma'$ such that $\mathfrak{X}'$ is a regular, strongly semi-stable arithmetic surface.
Such a $\Sigma'$ can be constructed as follows.
Let $\sigma_1, \dots, \sigma_k \in \Sigma(\tilde{\Delta})$ be the two-dimensional cones that are strictly contained
in the open upper half-space --- these correspond to the edges of $\tilde{\Delta}$ that are not projected on the boundary of $\Delta$.
Let $\Sigma_0 = \Sigma(\tilde{\Delta})$ and repeat the following for $i=1, \dots, k$.

\begin{itemize}
\item[] Because $v(\Delta \cap \mathbb{Z}^2) \subset \mathbb{Z}$, the extremal rays of $\sigma_i$ are generated by
vectors
\[ (\alpha, \beta, 1) \quad \text{and} \quad (\gamma,\delta,1) \quad \text{with $\alpha,\beta,\gamma,\delta \in \mathbb{Z}$},\]
hence by applying a $\mathbb{Z}$-affine transformation if necessary, we may
assume that $\sigma_i$ is generated by
\[ (0,0,1) \quad \text{and} \quad (d,0,1), \]
where $d$ is the greatest common divisor of the $(2 \times 2)$-minors
of
\[ \begin{pmatrix} \alpha & \beta & 1 \\ \gamma & \delta & 1 \\ \end{pmatrix}. \]
Subdivide $\sigma_i$ by introducing $d-1$ new rays, generated by the vectors 
$(1,0,1)$,$(2,0,1)$,$\dots$,$(d-1,0,1)$,
\begin{figure}[H]
\centering
\includegraphics[height=1.5cm]{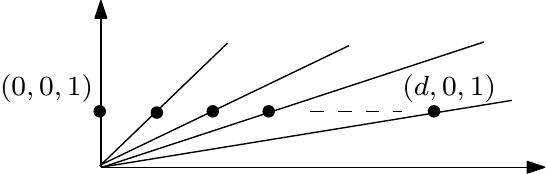}
\end{figure}
and denote the resulting fan by $\sigma_i^\text{sub}$.
Extend this to a subdivision of $\Sigma_{i-1}$ by connecting each of
the newly introduced rays with a fixed third extremal ray of each three-dimensional cone
adjacent to $\sigma_i$. Let $\Sigma_i$ be the resulting fan.
\begin{figure}[H]
\centering
\includegraphics[height=2.5cm]{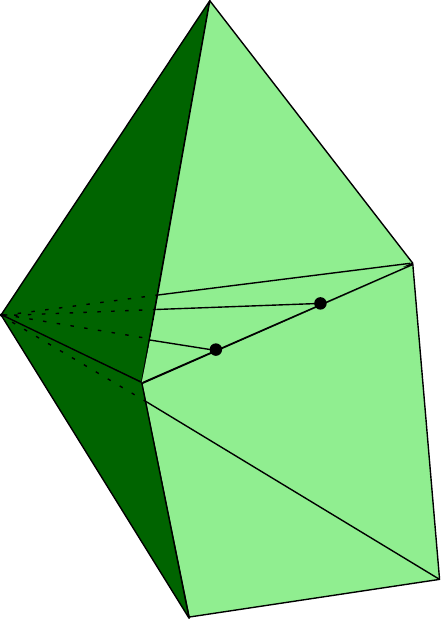}
\quad \quad $\stackrel{\longrightarrow}{\phantom{\begin{pmatrix} a \\ b \\ c \\ d \\ \end{pmatrix}}}$ \quad \quad
\includegraphics[height=2.5cm]{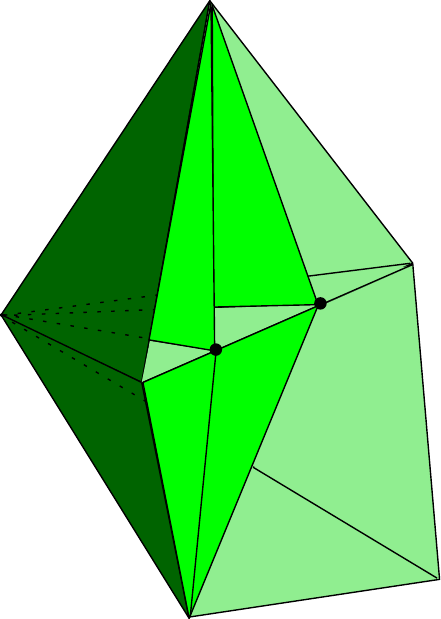}
\end{figure}
\end{itemize}
Then let $\Sigma' = \Sigma_k$. Note that the construction of $\Sigma'$ is highly non-canonical: it
depends on the way we ordered $\sigma_1, \dots, \sigma_k$ and
on the respective choices of the third extremal rays of the adjacent cones. But since by local non-degeneracy $X_0$ does not contain
any of the zero-dimensional toric orbits of $Y$, these choices affect neither $X'$ nor $p'$.

By normality, $Y$ is non-singular except possibly
at its one-dimensional and zero-dimensional toric orbits.
If $\tilde{\tau}$ is the graph of $v$ restricted to an edge $\tau$ of $\Delta$, then $Y$ is non-singular at $O(\tilde{\tau})$:
the corresponding cone of $\Sigma(\tilde{\Delta})$ is generated by vectors of the form $(a,b,0), (\alpha, \beta, 1)$ with $\gcd(a,b) = 1$, hence it is smooth.
By local non-degeneracy, we conclude that $X$ cannot have any singularities at $X_0$, except possibly at the toric
orbits $O(\tilde{\tau})$ associated to the lower edges $\tilde{\tau}$ of $\tilde{\Delta}$ that are not of
the above form $\text{graph}(v|_\tau)$. These edges exactly correspond to the cones $\sigma_1, \dots, \sigma_k$.

To prove that $\mathfrak{X}'$ is a regular strongly semi-stable arithmetic surface,
it suffices to make a local analysis around these toric orbits.
That is, for $i=1, \dots, k$ we consider the strict transform of $X \cap \text{Tor}(\sigma_i)$ under the restriction
of $\rho$ to $\text{Tor}(\sigma_i^\text{sub})$.
Let $\tilde{\tau}_i$ be the lower edge of $\tilde{\Delta}$ corresponding to $\sigma_i$.
As mentioned above, modulo a $\mathbb{Z}$-affine
transformation we may assume that $\sigma_i$ is generated by $(0,0,1)$ and $(d,0,1)$. In fact,
we can make the slightly stronger assumption that $\tilde{\tau}_i$ is supported on the $y$-axis, that the supporting hyperplanes of
the adjacent facets $\tilde{\Delta}_\ell$
and $\tilde{\Delta}_m$ contain $(1,0,0)$ resp.\ $(-1,0,d)$, and that the $t$-direction remains vertical. By local non-degeneracy, we can write
\[ f_t = g_{\tilde{\tau}}(y^{\pm 1}) + u \cdot g_{\tilde{\Delta}_\ell}(y^{\pm 1}, u) + x \cdot g_{\tilde{\Delta}_m}(y^{\pm 1}, x) +
t \cdot g_{\tilde{\Delta}}(y^{\pm 1},u,x,t), \quad u = x^{-1}t^d,\]
where
\begin{itemize}
\item $g_{\tilde{\tau}} \in \mathbb{C}[y^{\pm 1}]$ is a square-free Laurent polynomial (having $L(\ell,m)$ zeroes in $O(\tilde{\tau})$),
\item $g_{\tilde{\tau}} + u \cdot g_{\tilde{\Delta}_\ell} \in \mathbb{C}[y^{\pm 1},u]$
defines a smooth curve in $\mathbb{T}^2 = O(\tilde{\Delta}_\ell)$ (the completion
inside $\text{Tor}(\tilde{\Delta}_{\ell})$ of which is exactly $X^{(\ell)}$),
\item $g_{\tilde{\tau}} + x \cdot g_{\tilde{\Delta}_m} \in \mathbb{C}[y^{\pm 1},x]$ defines a smooth curve in $\mathbb{T}^2 = O(\tilde{\Delta}_m)$ (the completion inside $\text{Tor}(\tilde{\Delta}_m)$ of which is exactly $X^{(m)}$).
\end{itemize}
Then locally, $X$ is defined by $f_t(y^{\pm 1},u,t,x)$ inside
\[ \text{Tor}(\sigma_i) = \text{Spec} \ \frac{\mathbb{C}[y^{\pm 1},u,t,x]}{(t^d - ux)} \ \subset \text{Tor}(\Sigma(\tilde{\Delta})). \]
We will restrict our analysis of its strict transform in $\text{Tor}(\sigma_i^\text{sub})$ to the patch $\text{Tor}(\sigma_{i,1}^\text{sub})$, where
$\sigma_{i,1}^\text{sub}$ is the cone spanned by $(0,0,1)$ and $(1,0,1)$.
The dual cone
is generated by $(-1,0,1)$ and $(1,0,0)$, hence
\[ \text{Tor}(\sigma_{i,1}^\text{sub}) = \text{Spec} \ \mathbb{C}[y^{\pm 1},v,x] , \quad v = x^{-1}t.\]
  \begin{figure}[H]
  \centering
  \includegraphics[height=3cm]{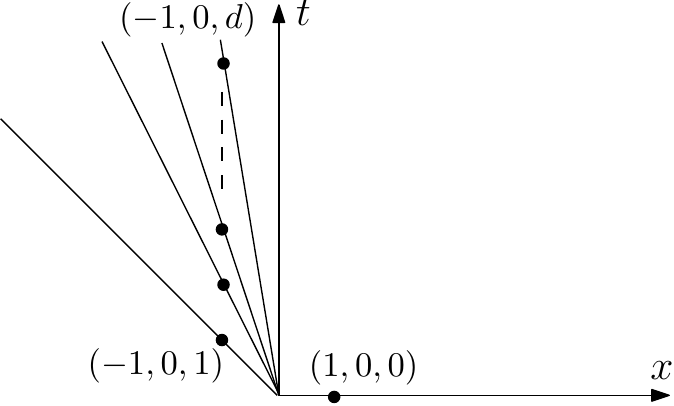}
  \end{figure}
\noindent Since this dual cone contains the dual cone of $\sigma_i$, we have a natural inclusion map which exactly describes
our toric resolution $\rho$:
\begin{equation} \label{resol}
  \text{Tor}(\sigma_{i,1}^\text{sub}) \rightarrow \text{Tor}(\sigma_i) : (y, v, x) \mapsto (y, v^dx^{d-1} , vx, x).
\end{equation}
The strict transform of $X$ under this map is described by
\begin{eqnarray*} f'_t(y^{\pm 1}, v, x) &=& g_{\tilde{\tau}}(y^{\pm 1}) + v^dx^{d-1} \cdot g_{\tilde{\Delta}_\ell}(y^{\pm 1}, v^dx^{d-1}) \\
&&+ x \cdot g_{\tilde{\Delta}_m}(y^{\pm 1}, x) + vx \cdot g_{\tilde{\Delta}}(y^{\pm 1},v^dx^{d-1},x,vx).
\end{eqnarray*}
The fiber above $t=0$ corresponds to taking $v=0$, in which case we find $g_{\tilde{\tau}}(y^{\pm 1}) +  x \cdot g_{\tilde{\Delta}_m}(y^{\pm 1}, x) = 0$
(the curve $X^{(m)}$), and taking $x=0$, in which case we find $g_{\tilde{\tau}}(y^{\pm 1}) = 0$ ($L(\ell,m)$ exceptional lines). These are easily
checked to be non-singular points of the strict transform, and all components intersect each other transversally. By making a similar analysis
of the other patches, one concludes that $\mathfrak{X}'$ is indeed a regular, strongly semi-stable arithmetic surface.

Now the generic fibers of $\mathfrak{X}$ and $\mathfrak{X}'$ are isomorphic, because $\rho|_{X'}$ is an isomorphism
on $p'^{-1}(V)$ for an open subset $V$ of $\mathbb{P}^1$. On the other hand, the special fiber of
$\mathfrak{X}$' differs from the special fiber of $\mathfrak{X}$. To see how the latter modifies under the above toric resolution, it
suffices to have a second look at the above analysis. Suppose that $\tau_i$ corresponds to
adjacent lower facets $\tilde{\Delta}_\ell$ and $\tilde{\Delta}_m$ of $\tilde{\Delta}$.
Then $d(\ell,m) - 1$ new rays are introduced. The introduction of a first new
ray separates the curves $X^{(\ell)}$ and $X^{(m)}$, and each intersection point becomes replaced
by an exceptional curve intersecting $X^{(\ell)}$ and $X^{(m)}$ transversally. This exceptional curve is contained
in the strict transform of $X$ and hence belongs to the special fiber of our new arithmetic surface. All intersections
remain transversal. More generally, if $d(\ell,m)-1$ rays are added, then each intersection point becomes replaced by a chain of
$d(\ell,m)-1$ transversally intersecting exceptional curves. Hence in the dual graph of $\mathfrak{X}$, if an edge corresponds to
$\Delta_\ell$ and $\Delta_m$, then it becomes replaced by a linear graph of length $d(\ell,m)$.

\section{Further remarks and adjustments}

\begin{itemize}
  \item In the third paragraph of the proof of \cite[Theorem~10]{polgon}, one should let $\mathfrak{X}'$ be the regular strongly semi-stable
  arithmetic surface constructed in Section~\ref{proof} above.

  \item The graph associated to the example following \cite[Conjecture~2]{polgon} should be replaced, e.g.\ by
  \begin{figure}[H]
  \centering
  \includegraphics[height=3cm]{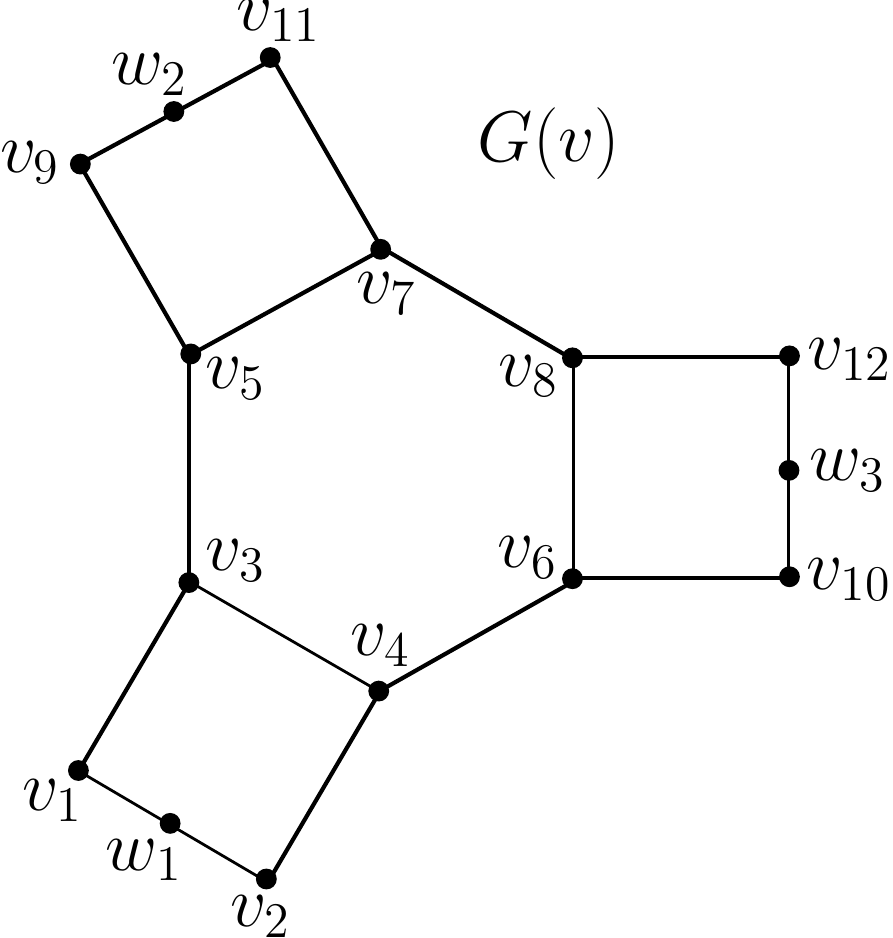}
  \end{figure}
  (we leave the determination of $v$ as an easy exercise) --- the according conclusion remains unaffected.

  \item In the proof of \cite[Theorem~11]{polgon}, it is easy to find an upper-convex piece-wise linear function $v$ realizing the given subdivision, such that $v(\Delta \cap \mathbb{Z}^2) \subset \mathbb{Z}$ and all $d(\ell,m)$'s equal $1$. Therefore, the proof remains valid.
\end{itemize}

\end{document}